\documentclass[a4paper,twoside, 10pt,final,leqno]{article} 
\usepackage[english]{babel}
\usepackage{amsthm,amsmath,amssymb,amscd,amsxtra,enumerate,amsgen}
\usepackage[all]{xy}
\usepackage[notref,notcite]{showkeys} 

\newcommand{\Ok}[1]{\ensuremath{{\mathcal O}_{#1}}}
\newcommand{\D}[1]{\ensuremath{{\mathcal D}_{#1}}}

\newcommand{\nc}{\newcommand}
\renewcommand{\frak}{\mathfrak}
\providecommand{\cal}{\mathcal}
\renewcommand{\bold}{\mathbf}

\numberwithin{equation}{section}

\newcommand{\pfname}{Proof.}

\newenvironment{pf}{\vskip-\lastskip\vskip\medskipamount{\it\pfname}}%
                      {$\square$\vskip\medskipamount\par}

\newenvironment{pfof}[1]{\vskip-\lastskip\vskip\medskipamount{\it
    Proof of #1.}}%
                      {$\square$\vskip\medskipamount\par}

\newtheorem{thm}{Theorem}[subsection]

\newtheorem{cor}[thm]{Corollary}
\newtheorem{corollary}[thm]{Corollary}

\newtheorem{prop}[thm]{Proposition}
\newtheorem{proposition}[thm]{Proposition}

\newtheorem{lemma}[thm]{Lemma} 

\theoremstyle{definition}
\newtheorem{defn}[thm]{Definition}
\newtheorem{definition}[thm]{Definition}

\theoremstyle{remark}

\newtheorem{remark}[thm]{Remark}

\nc{\Theorem}[1]{Theorem~{#1}}
\nc{\Th}[1]{({\sl Th.}~#1)}
\nc{\Theorems}[2]{Theorems~{#1} and ~{#2}}
\nc{\Thms}[2]{(\it Thms. ~{#1} and ~{#2})}

\nc{\manga}[6]{({\it Thms. ~ #1,\\ ~ #2, ~ #3, ~ #4, ~ #5, ~ #6})}

\nc{\Prop}[1]{({\sl Prop.}~{#1})}
\nc{\Proposition}[1]{Proposition~{#1}}
\nc{\Propositions}[2]{Propositions~{#1} and ~{#2}}
\nc{\Props}[2]{({\sl Props.}~{#1} and ~{#1})}
\nc{\Cor}[1]{({\sl Cor.}~{#1})}
\nc{\Corollary}[1]{Corollary~{#1}}
\nc{\Corollaries}[2]{Corollaries~{#1} and ~{#2}}
\nc{\Definition}[1]{Definition~{#1}}
\nc{\Defn}[1]{({\sl Def.}~{#1})}
\nc{\Lemma}[1]{Lemma~{#1}} 
\nc{\Lem}[1]{({\sl Lem.} ~{#1})} 
\nc{\Eq}[1]{equation~({#1})}
\nc{\Equation}[1]{Equation~({#1})}
\nc{\Section}[1]{Section~{#1}}
\nc{\Sec}[1]{({\sl Sec.} ~{#1})} 
\nc{\Chapter}[1]{Chapter~{#1}}
\nc{\Chapt}[1]{({\sl Ch.}~{#1})}

\nc{\Ex}[1]{{\sl Ex.}~{#1}}
\nc{\Exa}[1]{{\sl Example}~{#1}}
\nc{\Example}[1]{{\sl Example}~{#1}}
\nc{\Examples}[1]{{\sl Examples}~{#1}}
\nc{\Exercise}[1]{{\sl Exercise}~{#1}}

\nc{\Rem}[1]{({\sl Rem.}~{#1})}
\nc{\Remark}[1]{{\sl Remark}~{#1}}
\nc{\Note}[1]{{\sl Note}~{#1}}

\nc{\Conjecture}[1]{Conjecture~{#1}}
\nc{\Claim}[1]{Claim~{#1}}

\nc \Proof{{  \it Proof. }}

\nc{\xmu}{\mu}
\nc{\w}{\omega}

\nc \Ab{{\bold A}}
\nc \Gb{{\bold G}}
\nc \Qb{{\bold Q}}
\nc \Rb{{\bold R}} \nc \Cb{{\bold C}} 
\nc \Fb{{\bold F}}
\nc \Pb{{\bold P}}
\nc \SPb{{\bold {SP}}}
\nc \Zb{{\bold Z}} 
\nc \Nb{{\bold N}} 
\nc \Hb{{\bold H}} 
\nc{\wb}{\mathbf{w}}
\nc \syz{{\mathbf {syz}}}
\nc \mf{\frak m} \nc \mh{\hat{\m}} 
\nc \nf{\frak n}
\nc \Of{\frak O}
\nc \rf{\frak r}
\nc \mufr{{\mathbf \mu}}
\nc \hf{\frak h} 
\nc \qf{\frak q} 
\nc \bfr{\frak b} 
\nc \kfr{\frak k} 
\nc \pfr{\frak p} 
\nc \af{\frak a }
\nc \cf{\frak c }
\nc \sfr{\frak s} 
\nc \ufr{\frak u} 
\nc \g{\frak g} 
\nc \gA{\g_{\Ao}} 
\nc \lfr{\frak l}
\nc \afr{\frak a}
\nc \gfh{\hat {\frak g}}
\nc \gl{\frak { gl }}
\nc \Sl{\frak {sl}}
\nc \SU{\frak {SU}}
\nc{\Homf}{\frak{Hom}}

\newcommand{\on}{\operatorname}

\nc\row{\on {row\ }}
\nc\col{\on {col\ }}
\nc\nullo{\on {null\ }}
\nc \Ann {\on {Ann }}
\nc \Ass {\on {Ass \ }}
\nc \Coker {\on {Coker}}
\nc \Co{\on C}
\nc \Homo{\on {Hom}}
\nc \Ker {\on {Ker}}
\nc \No {\on N}
\nc \NN {\on {NN}}
\nc \NGo {\on {NG}}
\nc \Oo {\on O}
\nc \ch {\on {ch}}
\nc \rko {\on {rk}}
\nc \Sing {\on {Sing}}
\nc \Reg {\on {Reg}}
\nc \CoI {\on {CI}}
\nc \CoM {\on {CM}}
\nc \Gor {\on {Gor}}
\nc \Type {\on {Type}}
\nc \can {\on {can}}
\nc \Top {\on {T}}
\nc \Tr {\on {Tr}}
\nc \rel {\on {rel}}
\nc \tr {\on {tr}}
\nc \trdeg {\on {tr.deg}}
\nc \codim {\on {codim }}
\nc \coht {\on {coht}}
\nc \coh {\on {coh}}
\nc \embdim{\on {emb \ dim \ }}
\nc \qcoh {\on {qcoh}}
\nc \grade {\on {grade}}
\nc \hto {\on {ht}}
\nc \depth {\on {depth}}
\nc \prof {\on {prof}}
\nc \reso{\on {res}}
\nc \ind{\on {ind}}
\nc \prodo{\on {prod}}
\nc \coind{\on {coind}}
\nc \Con{\on {Con}}
\nc \Crit{\on {Crit}}
\nc \Der{\on {Der}}
\nc \Char{\on {Char}}
\nc \Ch{\on {Ch}}

\nc \Ext{\on {Ext}}
\nc \Eo{\on {E}}
\nc \End{\on {End}}
\nc \ad{\on {ad}}
\nc \Ad{\on {Ad}}
\nc \gr{\on {gr}}
\nc \Fo{\on {F}}
\nc \Gr{\on {Gr}}
\nc \Go{\on {G}}
\nc \GFo{\on {GF}}
\nc \Glo{\on {Gl}}
\nc \Ho{\on {H}}
\nc \CMo{\on {\CM}}
\nc \SCM{\on {SCM}}
\nc \hol{\on {hol}}
\nc{\sgd}{\on{sgd}}
\nc \supp{\on {supp}}
\nc \ssupp{\on {s-supp}}
\nc \singsupp{\on {singsupp}}
\nc \msupp{\on {msupp}}
\nc \spec{\on {spec}}
\nc \spano{\on {span \ }}
\nc \Max{\on {Max}}
\nc \Min{\on {Min}}
\nc \Mod{\on {Mod}}
\nc \Rad {\on {Rad}}
\nc \rad {\on {rad}}
\nc \rank {\on {rank}}
\nc \range {\on {range}}
\nc \Slo{\on {SL}}
\nc \soc {\on {soc}}
\nc \Irr {\on {Irr}}
\nc \Imo {\on {Im}}
\nc \SSo{\on {SS}}
\nc \lub{\on {lub}}
\nc \gldim{\on {gl.d.}}
\nc \pdo{\on {p.d.}} 
\nc \ido{\on {i.d.}} 
\nc \dSSo{\dot {\SSo}}
\nc \So{\on S}
\nc \Io{\on I}
\nc \Jo{\on J}
\nc \jo{\on j}
\nc \Ko{\on K}
\nc \PBW{\Ac_{PBW}}
\nc \Ro{\on R}
\nc \To{\on T}
\nc \Ao{\on A}

\nc \Do{{\on D}}
\nc \Bo{\on B}
\nc \Po{\on P}
\nc \Qo{\on Q}
\nc \Zo{\on Z}
\nc \U{\on U}
\nc \wt{\on {wt}}
\nc \Uh{\hat {\U}}
\nc \T{\on T}
\nc \Lo{\on L}
\nc{\dop}{\on d}
\nc{\eo}{\on e}
\nc{\ado}{\on{ad}}
\nc{\Tot}{\on{Tot}}
\nc{\Aut}{\on{Aut}}

\nc{\overrightleftarrows}[2]{\overset{#1}{\underset{#2}{\rightleftarrows}}}

\nc{\CCF}{\cal{CF}}
\nc{\CDF}{\cal{DF}}
\nc{\CHC}{\check{\cal C}}

\nc{\Cone}{\on{Cone}}
\nc{\dec}{\on{dec}}
\nc{\Diff}{\on{Diff}}
\nc{\dirlim}{\underset{\to}{\on{lim}}}
\nc{\dpar}{\partial}
\nc{\GL}{\on{GL}}
\nc{\CGr}{\cal{G}r}
\nc{\pr}{\on{pr}}
\nc{\semid}{|\!\!\!\times}
\nc{\Hom}{\on{Hom}}
\nc \RHom{\on {RHom}}

\nc \Proj{\mathrm {Proj\ }}

\nc{\Id}{\on{Id}}
\nc{\id}{\on{id}}
\nc{\Ima}{\on{Im}}
\nc{\invtimes}{\underset{\gets}{\otimes}}
\nc{\invlim}{\underset{\gets}{\on{lim}}}
\nc{\Lie}{\on{Lie}}
\nc{\re}{\on{Re }}
\nc{\Pic}{\on{Pic }}
\nc{\LPic}{\on{LPic }}
\nc{\Sch}{\on{Sch}}
\nc{\Sh}{\on{Sh}}
\nc{\Set}{\on{Set}}
\nc{\spo}{\on{sp\  }}
\nc{\Spec}{\on{Spec}}
\nc{\mSpec}{\on{mSpec}}
\nc{\Specb}{\bold {Spec}}
\nc{\Projb}{\bold {Proj}}
\nc{\Specan}{\on{Specan}}
\nc{\Spo}{\on{Sp}}
\nc{\Spf}{\on{Spf}}
\nc{\sym}{\on{sym}}
\nc{\symm}{\on{symm}}
\nc{\rop}{\on{r}}
\nc{\Td}{\on{Td}}
\nc{\Tor}{\on{Tor}}

\nc{\Artin}{\cal{A}rtin}
\nc{\Dgcoalg}{\cal{D}gcoalg}
\nc{\Dglie}{\cal{D}glie}
\nc{\Ens}{\cal{E}ns}
\nc{\Fsch}{\cal{F}sch}
\nc{\Groupoids}{\cal{G}roupoids}
\nc{\Holie}{\cal{H}olie}
\nc{\Mor}{\cal{M}or}

\nc{\CF}{\mathcal{F}}
\nc \Kc{{\cal K}}
\nc \Lc{{\cal L} }
\nc \lc{{\cal l}} 
\nc \CC{{\cal C}} 
\nc \Cc{{\cal C} }
\nc \Pc{{\cal P} }
\nc \Dc{\ensuremath{\mathcal D}}
\nc \Ac{{\cal A}} 
\nc \Bc{{\cal B}}
\nc \Ec{{\cal E}}
\nc \Fc{{\cal F}}

\nc \Mcc{{\cal M}} \nc \hM{\hat{\Mcc}} \nc \bM{\bar {\Mcc}} \nc
\hbM{\hat{\bar \Mcc}}  
\nc \Nc{{\cal N}}
\nc \Hc{{\cal H}} 
\nc \Ic{{\cal I}} 
\nc \Oc{\ensuremath{\mathcal O}}
\nc \Och{\hat{\cal O}} 
\nc \Sc{{\cal S} }
\nc \Tc{{\cal T}} 
\nc \Vc{{\cal V}} 
\nc{\CA}{{\cal A}}
\nc{\CB}{{\cal B}}
\nc{\C}{{\cal F}}
\nc{\Gc}{{\cal G}}
\nc{\CH}{{\mathcal H}}
\nc{\CI}{{\cal I}}
\nc{\CM}{{\cal M}}
\nc{\CN}{{\cal N}}
\nc{\CO}{{\cal O}}
\nc{\Rc}{{\cal R}}
\nc{\CT}{{\mathcal T}}
\nc{\CU}{{\cal U}}
\nc{\CV}{{\cal V}}
\nc{\CZ}{{\cal Z}}
\nc{\Homc}{{\cal {Hom}}}

\nc{\fa}{\frak{a}}
\nc{\fA}{\frak{A}}
\nc{\fg}{\frak{g}}
\nc{\fh}{\frak{h}}
\nc{\fI}{\frak{I}}
\nc{\fK}{\frak{K}}
\nc{\fm}{\frak{m}}
\nc{\fP}{\frak{P}}
\nc{\fS}{\frak{S}}
\nc{\ft}{\frak{t}}
\nc{\fX}{\frak{X}}
\nc{\fY}{\frak{Y}}

\nc{\bF}{\bar{F}}
\nc{\bCP}{\bar{\cal{P}}}
\nc{\bm}{\mbox{\bf{m}}}
\nc{\bT}{\mbox{\bf{T}}}
\nc{\hB}{\hat{B}}
\nc{\hC}{\hat{C}}
\nc{\hP}{\hat{P}}
\nc{\htest}{\hat P}

\nc{\nen}{\newenvironment}
\nc{\ol}{\overline}
\nc{\ul}{\underline}
\nc{\ra}{\to}
\nc{\lla}{\longleftarrow}
\nc{\lra}{\longrightarrow}
\nc{\Lra}{\Longrightarrow}
\nc{\Lla}{\Longleftarrow}
\nc{\Llra}{\Longleftrightarrow}
\nc{\hra}{\hookrightarrow}
\nc{\iso}{\overset{\sim}{\lra}}
\nc{\dsize}{\displaystyle}
\nc{\sst}{\scriptstyle}
\nc{\tsize}{\textstyle}
\nen{exa}[1]{\label{#1}{\bf Example.\ } }{}
\nen{rem}[1]{\label{#1}{\em Remark.\ } }{}
\nen{exer}[1]{\label{#1}{\em Exercise.\ } }{}

\begin{document}

\title{Geometric interplay between function subspaces and their rings
  of differential operators} \author{by\\[6pt] Rikard
  B{\"o}gvad\footnote{Department of Mathematics, Stockholm University,
    S-106 91, Sweden, email:rikard@matematik.su.se}\quad and Rolf
  K{\"a}llstr{\"o}m\footnote{Department of Mathematics, University of
    G{\"a}vle, S-801 76, Sweden, email:rkm@hig.se} }

\maketitle

\begin{abstract}
  We study, in the setting of algebraic varieties, finite-dimensional
  spaces of functions $V$ that are invariant under a ring $\Dc^V$ of
  differential operators, and give conditions under which $\Dc^V$ acts
  irreducibly.  We show how this problem, originally formulated in
  physics \cite{kamran-Milson-Olver:invariant,turbiner:bochner}, is
  related to the study of principal parts bundles and Weierstrass
  points \cite{EGA4,laksov-thorup}, including a detailed study of
  Taylor expansions. Under some conditions it is possible to obtain
  $V$ and $\Dc^V$ as global sections of a line bundle and its ring of
  differential operators. We show that several of the published
  examples of $\Dc^V$ are of this type, and that there are many more
  --- in particular arising from toric varieties.
\end{abstract}

\section{Introduction}
Assume that $V$ is a finite-dimensional vector subspace of a ring of
functions $A$. The problem to describe properties of the ring of
differential operators $\Dc_A^V$ on $A$ that preserve $V$, which is
studied in this paper, has its origin in quantum mechanics.  Such
invariant subspaces form the foundation of the theory of quasi-exactly
solvable quantum models, as formulated by Turbiner, Shifman,
Ushveridze et al. (see the survey \cite{olver:survey} and the
introduction to \cite{kamran-Milson-Olver:invariant}, as well as the
references given there).

It is easy to see for the polynomial algebra $A=\Cb[x_1, \dots ,
x_n]$, that $\Dc_A^V$ always acts irreducibly on $V$ (see
\cite{post-turbiner} for the one-variable case), and a similar result
was proved using different methods for real-analytic functions on a
domain of some $\Rb^n$ by Kamran-Milson-Olver in
\cite{kamran-Milson-Olver:invariant}. The latter authors introduce a
condition on $V$--- regularity --- that ensures that $\Dc_A^V$ acts
irreducibly on $V$.

There is also a Lie-algebraic approach, by directly finding
realizations of a Lie algebra as differential operators on a
(polynomial) ring preserving a vector space
\cite{kamran-Milson-Olver:invariant,turbiner:bochner}. It was noted in
\cite{gonzales-hurtubise-kamran-olver:quantification} that some of
these examples stem from a global variety with a line bundle. An
acquaintance with \cite{beilinson-bernstein,borel:Dmod} makes it
evident that many examples actually are local descriptions of rings of
global differential operators of line bundles on homogeneous spaces or
toric varieties.

Thus there is an interesting global slant to the problem of finding
invariant subspaces.  To understand it and also to extend the results
of \cite{kamran-Milson-Olver:invariant} to other spaces than open
subsets of ${\mathbf R}^n$, it is natural to study invariant subspaces
in the setting of sheaves of functions. Then the same global object
will unite many different embeddings of $V$ in different rings of
functions. The working out of the details of this more powerful and 
flexible way of approaching invariant subspaces, and exemplifying it,
is the main contribution of the present paper.

In \cite{EGA4} differential operators are defined using bundles of
principal parts. This point of view greatly simplifies the analysis of
regularity. Our approach has the added advantage that it makes the
results on points of Weierstrass subvarieties (points of inflexion and
osculating planes) by Laksov-Thorup \cite{laksov-thorup} available for the study of  differential operators.
Rather to our surprise it turns out that some of the ideas and results
are the same; for example the concept of regularity from
\cite{kamran-Milson-Olver:invariant} is also treated in
\cite{laksov-thorup}.  It will moreover be clear that on a nonsingular
variety in characteristic $0$ the two natural set-ups, one being based
on the sheaf of principal parts, the other on differential operators,
are in perfect correspondence.

We will now describe the content of the paper in more detail.  Central
to our whole approach is the {\it Taylor map} $\dop^n_V: \Oc_{X}\otimes_k V\to
\Pc^n_{X/k}(M)$ to the bundle of principal parts of a locally free
$\Oc_{X}$-module $M$ (this will give matrix-valued differential
operators), where $V$ is a subspace of $\Gamma(X,M)$; this is studied in
Chapter 2. We are first interested whether there exists an integer
$n_{inj}(x)$ giving fibre-wise injectivity $\dop^n_{V,x}: k_x\otimes_kV\to
k_x\otimes_{\Oc_x}\Pc_{X/k}^n(M)$ when $n\geq n_{inj}(x)$, the {\it
  separability}\/ of $V$.  Using Krull's theorem on the intersection
of powers of an ideal we get simple criteria for separability in
\Theorem{\ref{prop-sep}} and \Corollary{\ref{cor-sep}}.
\Theorem{\ref{n-fin-sing}} contains conditions guaranteeing the
existence of $n_{inj}$, such that $n\geq n_{inj}$ implies injectivity of
the Taylor map at all points.  These results extend by different means
the conditions given in \cite{laksov-thorup}, inter alia from smooth
to geometrically integral varieties, as well as bypasses the
complicated proof referred to in \cite{kamran-Milson-Olver:invariant}.
Our analysis of separability generalises the latter authors' use of
regularity; it reaches a more fundamental level, since we work with
more general varieties.  The surjectivity of the Taylor map is also
useful: define $n_{surj}(n^1_{surj})$ if the Taylor map is injective
at all points (points of height 1) for $n\leq n_{surj}$ (respectively
$n^1_{surj}$); it is related to jet-ampleness of a line bundle.

Using the definition of differential operators of Grothendieck
$\Dc_{X/k} (M)=\varinjlim \Hom(\Pc^n_{X/k}(M);M)$, in Chapter 3 we
apply the earlier results on the Taylor map to describe the action of
differential operators on $V\subset \Gamma(X,M)$.  In
\Theorem{\ref{key-lemma}} the surjectivity of the evaluation map
\begin{displaymath} W^n_x: \Dc^n(M)_x \to Hom_k(V, M_x)
   \end{displaymath}
   for differential operators on $M$ of order $n$ is shown to be
   equivalent to the injectivity of the Taylor map $\dop^n_{V,x}$.  In
   particular, if $V$ is separated, $W^n_x$ is surjective for high $n$
   \Th{\ref{surjective-seq}}. As a consequence the ring of
   differential operators $\Dc_{X/k}^{V}(M)$ that preserve $V\subset \Gamma(X,M)$
   acts irreducibly if $X$ is affine \Th{\ref{global-dec}}. This
   extends the algebraic counterpart of the main result of \cite[Th.
   4.8]{kamran-Milson-Olver:invariant}, bypassing Hodge algebra, from
   Zariski open subsets of $\Rb^n$ to real quasi-projective varieties,
   noting they are affine.  The chapter begins by proving that if the
   Taylor map $\dop_V^n$ is injective for some $n$, then it is
   actually injective for $n\geq |V|-1$ \Prop{\ref{v-1}}, a result also
   contained in \cite{laksov-thorup}.  Weierstrass subschemes $W(V)$
   are defined using the semi-continuous function $x\mapsto n_{inj}(x)$
   (when $X/k$ is nonsingular); $W(V)$ was constructed in a different
   way in \cite{laksov-thorup}. Our construction also gives a
   decreasing filtration of $W(V)$, and a candidate for a definition
   of Weierstrass subschemes on any scheme locally of finite type.
   
   Up to this point our theory has been local. As mentioned above many
   of the interesting examples, in the sense that they give rise to
   physical models, arise from Lie algebras that somehow act as
   differential operators both on an affine variety and a
   finite-dimensional vector space.  In \cite[p.\ 
   316]{kamran-Milson-Olver:invariant} the authors ask for an
   explanation of the ``significant mystery [that] is the connection
   [of the algebraic approach of \cite{kamran-Milson-Olver:invariant}]
   with the Lie algebraic approach of \cite{turbiner:bochner} to
   quasi-exactly solvable modules''.  By inspection many of these
   examples are seen to be restrictions to an affine open subset $U$
   from a global variety $X$ equipped with a line-bundle $\Lc$ (modulo
   the annihilator $\Ann V$). In Chapters 5 and 6 we review
   homogeneous spaces and toric varieties, and by combining results in
   the literature with our setup it is immediate to see that both
   these classes of varieties produce subspaces $V\subset A:={\mathbf
     C[x_{1},\ldots,x_{n}]}$ such that $\Dc^V$ (${\operatorname {mod}} \Ann V$) arise
   as restrictions of the ring of global differential operators on an
   ample line bundle.  Hence the Lie algebras arise from an invisible
   global context, while the approach of
   \cite{kamran-Milson-Olver:invariant} is local. This explains, we
   believe the connection. It also makes it interesting to consider
   the question of when a vector space $V\subset \Gamma(X, \Oc_{X})$
   and its ring of differential operators $\Dc^V_X$ comes by
   restriction from a global context.  There is an algebraic geometric
   technique, known to the ancients, to use $V$ to define a projective
   variety $X\subset X_{V}$ and a line-bundle $\Lc$ on $X_{V}$; in
   Chapter 4 we describe it and give in \Theorem{\ref{ext-prop}}
   conditions when we may extend differential operators $\Dc_{X}^{V}$
   to differential operators on $X_{V}$ that preserve $\Lc$. This is
   used in \Proposition{\ref{hidden-symmetry}} to prove that this
   machinery will detect ``hidden symmetries'': if $V\subset A$ is a
   finite-dimensional vector space such that there is a reductive Lie
   subalgebra $\g$ of $\Dc_{A}$ that acts irreducibly on $V$ and
   locally transitively on $A$ with a parabolic sub-algebra
   stabilising a point, then $\g $ is the restriction of differential
   operators on $X_{V}$ that preserve $\Lc$.
   
   Musson \cite{musson} describes the ring of global differential
   operators on a line bundle.  In particular, it is not difficult to
   use this description to give many examples of smooth varieties
   where the ring of differential operators preserving a vector space
   $V$ is not generated by first order differential operators.  These
   general results and toric constructions $X_V$, corresponding to a
   $V$ with a basis of monomials, do not seem to have been used before
   for the study of $\Dc^V$, even though Hirzebruch surfaces have been
   studied
   \cite{gonzalez-lopez-kamran-olver:quantification,finkel-kamran}.

   We have further used both toric varieties and homogeneous spaces to
   exemplify our concepts, in particular we use the equivariant
   structure on the bundle of principal parts
   \Prop{\ref{G-linearized}} to determine (for toric varieties) and
   estimate (for homogeneous varieties) $n_{inj}$, and in the case of
   Hirzebruch surfaces, we also make more explicit point-wise
   calculations.
   
   We would like to thank T. Ekedahl for valuable comments.

   {\it Notations and assumptions}:\/ A variety $X/k$ is an
   irreducible scheme of finite type over a field, and $\Oc_X$ is its
   structure sheaf. If $x$ is a point on a scheme $X$ and $M$ a sheaf
   of modules over $\Oc_X$, we let $k_x$ be the residue field of
   $\Oc_X$, $M_x$ the stalk of $M$, and $k_x\otimes_{\Oc_x}M_x$ the fibre of
   $M$, at $x$. If $\phi : \Fc \to \Gc$ is a map of $\Oc_X$-modules, then
   $\phi_x : \Fc_x \to \Gc_x$ is the map of stalks and $\phi(x):
   k_x\otimes_{\Oc_x}\Fc_x \to k_x \otimes_{\Oc_x}\Gc_x$ the map of fibres at $x$.
   Throughout the paper the base field $k$ is of characteristic $0$
   and $M$ will denote a locally free $\Oc_X$-module of finite rank.

\section{Taylor maps}
\subsection{Generalities}\label{set-up}
A general reference for the material in this section is \cite{EGA4}.
Let $(X/k,\Oc_X)$ be a variety and $M$ a locally free $\Oc_X$-module
of finite rank.  Let $\Delta: X \to X\times_k X$, $x\mapsto (x,x)$, be the
diagonal map,
$I_\Delta$ be the kernel of the mapping $\Delta^*(\Oc_{X\times_k X})\to
\Oc_X$, and
put, for each integer $n\geq 0$, $\Pc_{X/k}^n= \Delta^*(\Oc_{X\times_k
  X})/I^{n+1}_\Delta$; define also $\Pc_{X/k}^{-1}=0$.  The sheaf of
$(\Oc_X, \Oc_X)$-bimodules $\Pc^n_{X/k}$ is the sheaf of $n$th order
principal parts.  Put also $\Pc^n_{X/k}(M)= \Pc^n_{X/k}\otimes_{\Oc_X}M$.

Let $V $ be a finite-dimensional $k$-subspace of the space of global
sections $M(X):=\Gamma(X,M)$; $V$ is also regarded as a constant sheaf on
$X$. There is a map $M \to \Pc^n_{X/k}(M)$, $m\mapsto 1\otimes m$, which is
injective since $M$ is locally free, and composing it with the
injective map $V \to M$ we get an injective map $ V \to \Pc^n_{X/k}(M)$.
Putting $\Vc_X=\Oc_X\otimes_k V$ we get a map $\dop_V^n :\Vc_X \to
\Pc^n_{X/k}(M) $, $\phi\otimes v \mapsto \phi \otimes v \mod
I_\Delta^{n+1}$; this map however need
not be injective.  Let $\Kc^n$ and $\Cc^n$ be the kernel and cokernel,
respectively, of $\dop_V^n$, so we have the exact sequence, which we
will refer to as the {\it Taylor sequence},
\begin{equation}\label{w-sequence}
   0\to \Kc^n \to \Vc_X \xrightarrow{\dop^n_V} \Pc^n_{X/k}(M) \to
\Cc^n\to 0.
\end{equation}
Varying the integer $n$ one gets different exact sequences
(\ref{w-sequence}), connected in the diagram
\begin{displaymath}
  \xymatrix{
0\ar[r]&\Kc^{n+1} \ar[dd]\ar[dr]&&\Pc^n_{X/k} \ar[r]& \Cc^n \ar[r]&0\\
&&\Vc_X\ar^{\dop^n_V}[ur]\ar^{\dop^{n+1}_V}[dr]&&&\\
0\ar[r]&\Kc^n \ar[ur]&&\Pc^{n+1}_{X/k} \ar[uu]_{q_{n+1}}\ar[r]&
\Cc^{n+1}\ar[uu] \ar[r]&0.\\
}
\end{displaymath}
where $q_{n+1}$ is the natural projection map.  The inverse limit
$\Pc^\infty_{X/k}(M):= \varprojlim_{ n } \Pc_{X/k}^n(M)$ is provided with
the $I_\Delta$-adic topology, which is used to define the sheaf of
differential operators as the $\Oc_X$-bimodule of continuous maps
\begin{displaymath}
\Dc_X(M):=Hom_{\Oc_X}^{cont}(\Pc^\infty_{X/k}(M),M).
\end{displaymath}
where $M$ is given the discrete topology.  The $\Oc_X$-bimodule of
differential operators of order at most $n$ is denoted $\Dc_X^n(M)$,
which together for all $n$ give a filtration
\begin{eqnarray*}
   \Dc_X^n(M)&:=&Hom_{\Oc_X}(\Pc_{X/k}^n(M),M), \\
0 \subset \Dc_X^0(M)= End_{\Oc_X}(M)&\subset& \cdots \subset
\Dc_X^n(M)\subset \Dc_X^{n+1}(M)\subset \cdots, \\
\Dc_X(M)&=& \cup_{n=0}^\infty\Dc^n_X(M).
\end{eqnarray*}
$M$ is $>1$.  The fibre of the sheaf of principal parts at a
$k$-rational point $x$ in $X$ is
\begin{equation}\label{princ-fibre}
   k\otimes_{\Oc_x}\Pc^n_{X/k}(M)_x \cong M_x/\mf^{n+1}_x M_x\quad
(\cite[16.4.11]{EGA4})
\end{equation}
and if the stalk $\Pc^n_{X/k}(M)_x$ is free, then the fibre
$k\otimes_{\Oc_x} \Dc_x^n(M) = (M_x/\mf^{n+1}_x M_x)^*:=
Hom_{k}(M_x/\mf^{n+1}_x M_x, k)$.

We will need that $\Pc_{X/k}^n(M)$ be locally free over $\Oc_X$ (of
finite rank).
\begin{lemma}\label{locfree}
  A scheme of finite type $X/k$ is regular if and only if the
  $\Oc_X$-module $\Pc^n_{X/k}$ is locally free for each integer $n$.
    \end{lemma}
\begin{pf}  If $X$ is regular, $X$ is smooth, and $\Omega_{X/k}$ is
  locally free (\cite[Cor.  17.5.2,17.15.6]{EGA4}). Hence
  $\Pc_{X/k}^n$ is locally free (\cite[Thm.  16.12.2,
  Def.16.10.1]{EGA4}). The converse follows from the same reference.
\end{pf}

Let
\begin{displaymath}
   \dop^n_{V,x}:\Vc_x \to \Pc_{X/k,x}^n(M)
\end{displaymath}
be the map of stalks at $x$. This is the $n$th Taylor expansion map at
$x$ of the vectors in $V$. In a smooth rational point $x$, there is a
basis of $\Pc_{X/k,x}^n(\Oc_X)$, consisting of all monomials of degree
less than $n$ in $\xi_{i}=x_{i}\otimes 1-1\otimes x_{i}$, where $x_{i}$ are a
regular system of parameters at $x$. Let the derivations $\partial_{1},\ldots,
\partial_{d}$ correspond to the coordinates; then the map $\dop^n_{V,x}$ is
described by the matrix with rows $(\partial^\alpha (m_i)/\alpha!)$, for
$m_i$ a basis
of $V$. Here the multi-indices $\alpha=(\alpha_{1},\ldots,\alpha_{d})$ take
all possible
values that define a differential operator of order less than or equal
to $n$, $\partial^\alpha=\prod(\partial/\partial x_{i})^{\alpha_{i}}$ and
$\alpha!=\prod \alpha_{i}!$.

The map of fibres at a point $x$ is with residue field $k_x$ is
\begin{displaymath}
\dop^n_V(x):k_x\otimes _kV \to k_x\otimes_{\Oc_x}\Pc^n_{X/k,x}(M).
\end{displaymath}

\subsection{The injectivity of the Taylor map}
The following elementary fact will play a central part in the paper
\cite[Proposition II.3.6]{bourbaki:commutative}.

\begin{prop}\label{splitting}
  Let $\phi :\Fc \to \Gc$ be a map of locally free $\Oc_X$-modules of
  finite rank. Then the following are equivalent at a point $x$ in
  $X$:
   \begin{enumerate}
   \item the map of fibres $\phi(x): k_x\otimes_{\Oc_x}\Fc_x \to
     k_x\otimes_{\Oc_x}\Gc_x$ is injective;
   \item the map of stalks $\phi_x :\Fc_x \to \Gc_x$ is split injective.
   \end{enumerate}
\end{prop}

We will employ Krull's theorem in order to get separated topologies on
a Noetherian tensor product $K_1\otimes_k A$.  We need therefore that
$K_1\otimes_k A$ be integral, so the following lemma is useful.
\begin{lemma}(\cite[Cor.
  4.6.3]{EGA4}) \label{integraltensor} Let $X/k$ be an integral
  scheme.  The following are equivalent:
  \begin{enumerate}
  \item $k$ is algebraically closed in the function field of $X$;
  \item $X\times_{k} K$ is integral for each extension $K$ of $k$.
\end{enumerate}
\end{lemma}
One says that the scheme $X/k$ is {\it geometrically integral}\/ if it
satisfies these conditions. In particular, $\Oc_x\otimes_{k}K$ is always
integral.

That the Taylor map is injective will turn out to have great
significance for the properties of differential operators that we are
interested in. We have an exact sequence
\begin{displaymath}
  0 \to \Kc \to \Delta^*(p^*_2(M)) \xrightarrow{\dop}
\Pc_{X/k}^{\infty}(M)
\end{displaymath}
where $\Kc= \cap_{n\geq 1} I^n_\Delta M$ and $p_2 :X\times_k X \to X$ is
the projection
on the second factor; $I^n_\Delta M$ denotes the image of the canonical map
$I^{n+1}_\Delta \otimes_{\Oc_X}M \to \Delta^*(p^*_2(M))$ (which is
injective when
$\Pc^n_{X/k}$ is flat over $\Oc_X$). Hence $\dop$ is injective if the
$I_\Delta$-adic topology on the $\Oc_X$-module $\Delta^*(\Oc_{X\times_k
  X}\otimes_{\Oc_X}M)$ is separated.
\begin{prop}\label{taylorinjective}
  Let $X/k$ be a scheme of finite type.  If $X\times_k X$ is integral (e.g.
  $X/k$ is geometrically integral), then the $I_\Delta$-adic topology on
  $\Delta^*(p_2^*(M))$ is separated, or equivalently, the Taylor map $\dop$
  is injective.
\end{prop}
\begin{pf}
  The map $\dop$ is injective if it is injective at associated points
  of $X$. It therefore suffices to prove the following.  Let $K/k$ be
  a field extension and $M$ be a finite-dimensional linear space over
  $K$ and $I_\Delta$ be the kernel of $K\otimes_k K \to K$, $a\otimes b \mapsto ab$. Then the
  map $d_{K/k,M}: K\otimes_k M \to \Pc^\infty _{K/k}(M)= \varprojlim_{n\to \infty} K\otimes
  K/I^{n+1}_\Delta \otimes_K M$ is injective.  But $\Ker(K\otimes_kK \to \Pc^\infty _{K/k})=
  \cap_{n\geq 1} I_K^n$ and $\Ker(K\otimes_kM \to \Pc^\infty_{K/k}(M))= \cap (I^n_\Delta M)=
  (\cap_{n\geq 1} I_K^n)M$ (recall $M$ is a linear space).  Since $K/k$ is
  of finite type it follows that $K\otimes_k K$ is noetherian; by assumption
  it is also integral, hence by Krull's theorem describing the
  intersection of powers of an ideal \cite[{Proposition
    III.3.5}]{bourbaki:commutative} $(\cap_{n\geq 1} I_K^n)M=0$ so
  $d_{K/k,M}$ is injective.
\end{pf}
A simple example where the Taylor map $\dop_V$ never is injective is
provided by any non-trivial finite algebraic extension $K/k$, with
$V=K=M$; since $\Pc^n_{K/k}=K$, the Taylor map will never be
injective.

\Proposition{\ref{taylorinjective}} implies easily that there exists
an integer $N$ such that if $n\geq N$, then the map
\begin{displaymath}
\dop^n_V : \Vc_X \to \Pc_{X/k}^n(M)
\end{displaymath}
is injective in this situation. We will refine this result by
investigating conditions that may be put on $V$ instead of $X/k$.

\begin{prop}\label{sep-ring-prop}
  Let $R$ be a local noetherian integral $k$-algebra such that its
  residue field $l$ is a finitely generated extension of $k$, $M$ a
  free $R$-module and $V$ finite-dimensional $k$-subspace of $M$. Let
  $k^r$ be the algebraic closure of $k$ in the quotient field $K=
  K(R)$ of $R$ and $k^l$ its algebraic closure in $l$. Let
  $R^r=k^rR\subset K$ be the $k^r$-algebra generated by $k^r$ and $R$,
  and $M^r=k^rM$ the $k^r$-space generated by $M$ in $K\otimes_R M$.
  The following are equivalent:
  \begin{enumerate}
  \item There is an integer $N$ such that
    $$
    l\otimes_{k}V\to l\otimes_R\Pc^n_{R/k}\otimes_R M
    $$
    is injective when $n\geq N$;

  \item The map
    \begin{displaymath}
      l\otimes_k V \to l\otimes_R\Pc^\infty_{R/k}\otimes_R M
    \end{displaymath}
    is injective;
  \item The map
  \begin{displaymath}
    k^l \otimes_k V \to \frac {k^l \otimes_k k^r}{J_0}
    \otimes_{k^r} M^r
  \end{displaymath}
  is injective, where
\begin{displaymath}
  J_0 = \{ y \in k^l \otimes_k k^r \subset k^l \otimes_k R^r \ \vert
\ \exists m \in I, \text{ s.t. } (1+m)y=0\}
\end{displaymath}
and
\begin{displaymath}
  I = \Ker (l\otimes_k R \to l).
\end{displaymath}
We have
\begin{displaymath}
  J_2 \subset J_0 \subset J_1,
\end{displaymath}
where $J_1 = \Ker (k^l \otimes_k k^r \to l\otimes_R R^r)$ and $J_2 $ is the
ideal in
$k^l \otimes_k k^r$ generated by $\Ker (k^l \otimes_k \hat k \to k^l)$,
where $\hat
k$ is the algebraic closure of $k$ in $R$.
  \end{enumerate}
\end{prop}
\begin{pf} It suffices to prove this when $M= R$.
  There is an injective map $l\otimes_R \Pc^\infty_{R/k}\to \varprojlim_n\
l\otimes_R
  \Pc^n_{R/k}$, so $(2)$ holds if and only if the map $l\otimes _k V \to
  \varprojlim_n \ l\otimes_R \Pc^n_{R/k}$ is injective. This is a
  reformulation of (1), since $V$ is finite-dimensional.  Furthermore
  we have a short exact sequence $0 \to I^n \to l\otimes_k R\to l\otimes_R
\Pc^n_{R/k} \to
  0$ (where $I=\Ker(l\otimes_k R\to l)$), and so $\Ker (l\otimes_k R \to
  \varprojlim_n\ l\otimes_R \Pc^n_{R/k}) = \cap I^n $. Put $J= \cap I^n$,
so $(2)$
  holds if and only if the map $l\otimes_k V \to \frac{l\otimes_k R}{J}$ is
  injective. Since $l/k$ is of finite type the ring $l\otimes_k R$ is
  noetherian by Hilbert's basis theorem and by Krull's theorem the
  ideal $J$ consists of the elements $s$ in $l\otimes_k R$ such that there
  exists $x\in I$ satisfying $(1+x)s=0$. Consider the inclusions
    \begin{displaymath}
      S:= l\otimes_k R \   \subset \ \bar S = l\otimes_k R^r \
\subset \ S_K = l\otimes_k K.
    \end{displaymath}
    We have $R \subset S$ and $R\cap J =0$, since elements in $R$ are not
    zero-divisors in $S$. Hence $\bar J: = \bar S J$ and $J_K: = S_K
    J$ are proper ideals.  Furthermore $S\cap J_{K}=J$; for if
$y/f=s\subset S\cap
    J_{K} $, there exists $1+x\in 1+I$, such that $(1+x)y=f(1+x)s=0$,
    hence $(1+x)s=0$ and so $s\in J$.  This means that there are maps
\begin{displaymath}
  S/J \to \bar S/ \bar J \to S_K/J_K,
\end{displaymath}
such that the composition is injective. Hence $(2)$ is equivalent to
the injectivity of
\begin{displaymath}
  l\otimes_k V \to \frac{l\otimes_k K}{J_{K}}.
\end{displaymath}
Now $ S_{K} =l\otimes_{k^l}(k^l\otimes_{k}k^r)\otimes_{k^r}K $. If $P$ is a
prime ideal
of $S_{0}:= k^l\otimes_k k^r$, then $S_K/S_KP=l\otimes_{\tilde
l}(S_{0}/P)\otimes_{k^r}K
$ is integral \Prop{\ref{integraltensor}}; hence $S_KP$ is a prime
ideal.  According to \cite[IV Prop 11]{bourbaki:commutative}, if $0=\cap
P_{i}$ in $S_{0}$, we then also have $0=\cap S_K P_{i}$ ($S_{0}$ and $S$
are reduced); this implies that the minimal primes of $S_{0}$ and
$S_K$ correspond, and more generally $A = S_K(A\cap S_{0})$ is true for
all radical ideals. Since all ideals in $S_K$ (a direct sum of fields)
are radical, this is true for all ideals.  Define $J_{0}:=S_{0}\cap
J_{K}$; then the equivalence of (2) with (3) is immediate:
$S_K/J_K=l\otimes_{k^l}(S_{0}/J_{0})\otimes_{k^r}K $, and the map $
l\otimes_k V \to
\frac{l\otimes_k K}{J_{K}}$, is a flat extension of the map $k^l\otimes_k V \to
(S_{0}/J_{0})\otimes_{k^r}K$. The characterisation of $J_{0}$ follows from
the fact that $J_{K}=\{ y \in S_{K} \ \vert \ \exists m \in I, \text{ s.t. }
(1+m)y=0\}$.

To estimate $J_{0}$, note first that $A:=\Ker (k^l \otimes_k \hat k \to
k^l)\subset
I$, that $A^n=A$, since $k^l \otimes_k \hat k $ is a direct sum of fields,
and that hence $J_{2}=S_{0}A\subset J_{K}$, implying that $J_{2}\subset J_{0}$.
Secondly put $\bar I =\bar S I= \Ker (l\otimes_k k^rR \to l\otimes_R R^r)$.
We have
an inclusion
\begin{displaymath}
  \bar J = \bar S J\  \subset
  \ B:= \cap \bar S I^n=  \ \cap \bar I^n
\end{displaymath}
Since $k^l \otimes_k k^r\subset l\otimes_k R^r=\bar S$, we can define $J_1
= \Ker (k^l
\otimes_k k^r \to l\otimes_R R^r)\subset \bar S$.  Then $J_{1}=S_{0}\cap
\bar I\subset \bar I$, and
hence, as before, $J_{1}\subset B=\cap \bar I^n$.  We assert that
$J_{K}\subset S_{K}
\cap B=S_{K}J_1$. The inclusion is immediate, and the description of $B$
follows, in the same way as above, since $J_{1}=S_{0}\cap (\cap \bar I^n)$,
and since $\bar S\cap S_{K}B=B$, by the characterisation $B=\{ y \in S_{K} \
\vert \ \exists m \in \bar I, \text{ s.t. } (1+m)y=0\}$ (cf.\ the similar
argument above)
\end{pf}

By reduction to $\Oc_x=R$ and $ M= M_x$, we get the following analysis
of the behaviour of the Taylor map at a point $x$.

\begin{thm}\label{prop-sep} Assume $X$ is an integral variety, $M$ a
  locally free $\Oc_X$-module and $V$ a $k$-subspace of $\Gamma(X,M)$, and
  let $x$ be a point in $X$. Denote the algebraic closure of $k$ in
  the function field $k(X)$ by $k^r$, and the $k^r$-module generated
  by $M_{x}$ in $k(X)\otimes_{\Oc_x} M_x$ by $M^r_x$. The following are
  equivalent:
\begin{enumerate}
\item There is an integer $N$ such that the map
  $$
  k_{x}\otimes_{k}V\to k_x\otimes_{\Oc_x}\Pc^n_{X/k}(M)_{x}
  $$
  is injective when $n\geq N$

\item The map
  $$
  k^l_{x}\otimes_{k}V\to ((k^l_{x}\otimes_{k}k^r)/J_{0}) \otimes_{ k^r} M_x^r
  $$
  is injective, where $k^l_{x}$ is the algebraic closure of $k$ in
  $k_{x}$ and $J_{0}$ is the ideal defined in the following way:
  $$J_{0}=\{y\in k^l_{x}\otimes_{k}k^r \ \vert \ \text{there is}\ m\in I \text{
    such that } (1+m)y=0\}.$$
\end{enumerate}
\end{thm}

\begin{definition} \label{sep-def}
  If the equivalent conditions in \Proposition{\ref{sep-ring-prop}}
  obtains we say that $V$ is a separated subspace of $M$.  If the
  equivalent conditions in \Theorem{\ref{prop-sep}} obtains at $x$
  we say that $V$ is separated at $x$, and if these conditions obtains
  at each point we say that $V$ is separated (on $X$).
\end{definition}

\begin{corollary} \label{cor-sep}
  Let $X/k, V$ and $M$ be as in \Theorem{\ref{prop-sep}}.
\begin{enumerate}
\item If $X$ is geometrically integral, then $V$ is separated at all
  points.
\item If $x$ is a rational point, then $V$ is separated at $x$.
\item If $x$ is a normal point, then $V$ is separated if and only if
  the map
  $$
  k^r\otimes_{k}V\to M_x^r
  $$
  is injective.
\item If the map $k^r \otimes_k V \to M_x^r$ is injective, then $V$ is
  separated at $x$.
        \end{enumerate}
\end{corollary}

\begin{pf}(1): If $X$ is geometrically integral then
  $k^r=k$ \Prop{\ref{integraltensor}} so $k^l_{x}\otimes_{k}k^r= k^l_{x}$
  and since $J_{0}$ is proper it is hence the zero-ideal. Thus the map
  in part (2) of \Theorem{\ref{prop-sep}} is just the inclusion
  $V\to M_{x}$, tensored by the flat $k$-module $k^l_{x}$.  (2): If $x$
  is a rational point $k^l_{x}=l$ and hence again $J_{0}=0$, and the
  map becomes the inclusion.  (3): If $x$ is normal then $\bar
  \Oc_{x}=\Oc_{x}$ and $\hat k=k^r$, and hence $J_{1}=J_{2}=J_{0}$,
  and so $\frac{k^l_{x}\otimes_{k}k^r}{J_{0}}=k^l_{x}$. This means that the
  sequence of the proposition becomes $k^l_{x}\otimes_{k}V\to k^l_{x}
\otimes_{k^r}
  k^r M_{x}$. This is by flatness injective if and only if the map
  $k^r\otimes_{k}V\to M^r$ is injective.  (4): From the assumption, it follows
  that the map
  $$
  \frac{k^l_{x}\otimes_{k}k^r}{J_{0}}\otimes_{k^r} k^r \otimes_k V \to
  \frac{k^l_{x}\otimes_{k}k^r}{J_{0}}\otimes_{k^r} M_x^r
  $$
  is injective. However, the first vector space is isomorphic to
  $\frac{k^l_{x}\otimes_{k}k^r}{J_{0}}\otimes_k V$, and contains
$k^l_{x}\otimes_{k} V$.
  Hence (4) follows from \Proposition{\ref{sep-ring-prop}}.
        \end{pf}

\begin{remark}\label{sep-ass}
  The following are equivalent:
  \begin{enumerate}
  \item $V$ is a separated subspace of $M_x$ at each associated point
    $x$;
  \item there exist an integer $N$ such that the Taylor map
    $\dop^n_{V}:\Vc \to \Pc^n_{X/k}$ is injective when $n\geq N$.
  \end{enumerate}
\end{remark}
\begin{definition}\label{def-inj}
  Let $n_{inj}(x)= n_{inj}(x,V)$ be the smallest integer such that
  $\dop^n_V(x)$ is injective (injectivity order of $V$ at $x$). Define
  also
  \begin{eqnarray*}
&n_{inj}= n_{inj}(V)= \sup \{n_{inj}(x)\ \vert \ x\in X \}, \\
&N_{inj}= N_{inj}(V)=\inf\{n\ \vert \ \dop^n_V \text{is injective}\}.
  \end{eqnarray*}
  The integer $n_{inj}$ is the injectivity order of $V$ and $N_{inj}$
  the generic injectivity order.
\end{definition}

\begin{prop}\label{semi-cont-reg} ($X/k$ is a smooth scheme locally
  of finite type) Let $M$ a locally free $\Oc_X$-module and $V$ a
  sub-space of $\Gamma(X,M)$. Then the function $x\mapsto n_{inj}(x)$ is upper
  semi-continuous.  In particular, if $x$ is a specialisation of a
  point $y$, then $n_{inj}(x)\geq n_{inj}(y)$, and if $X$ is irreducible,
  then for any point $x$ we have $n_{inj}(x)\geq N_{inj}$.
\end{prop}
We will see in \Proposition{\ref{v-1}} that $N_{inj}\leq \dim V-1$.
\begin{pf}
  Let $x$ be a specialisation of the point $y$, and let $n\geq
  n_{inj}(x)$. By \Proposition{\ref{splitting}}, $(1)\Rightarrow (2)$ the map
  $\dop^n_{V,x}$ is split injective, hence it is split injective in
  some neighbourhood $\Omega$ of $x$; since $y$ specialises to $x$ we have
  $y\in \Omega$, implying that the map $\dop^n_{V,y}$ is split injective,
  hence by \Proposition{\ref{splitting}}, $(2)\Rightarrow (1)$, $n\geq
n_{inj}(y)$.
  It is also clear that sets of the form $\{x\in X \ \vert \ n_{inj}(x)<
  i \}$ are open for each integer $i$. Since $\Vc_X$ is free, $N_{inj}
  = \sup \{n_{inj}(\eta)\ \vert \ \eta \text{ is an associated}$ $\text{ point of }
  X\}$, so in particular $n_{inj} (x)\geq N_{inj}$, $x\in X$, when there is
  only one associated point.
\end{pf}
By \Proposition{\ref{semi-cont-reg}} the function $x\mapsto n_{inj}(x)$ is
upper semi-continuous when $X$ is a regular variety, hence $n\geq
n_{inj}$ if $\dop^n_{V}(x): k_x \otimes_k V \to
k_x\otimes_{\Oc_x}\Pc^n_{X/k}(M)_x$
is injective at each closed point $x$. Actually, if $X/k$ is any
scheme locally of finite type with trivial Jacobson radical it follows
from the Nullstellensatz that $n\geq n_{inj}$ if $\dop^n_{V}(x): k_x \otimes_k
V \to k_x\otimes_{\Oc_x}\Pc^n_{X/k}(M)_x$ is injective at each closed point
$x$.

\begin{thm}\label{n-fin-sing}
  Let $ X/k$ be a reduced scheme locally of finite type, $M$ a locally
  free $\Oc_X$-module, and $V$ be a finite-dimensional separated
  $k$-sub-space of $\Gamma(X,M)$ (e.g.  $k$ is algebraically closed and $V$
  injects to each stalk $M_x$).  Then the injectivity order is finite,
  $n_{inj}(V)< \infty $.
\end{thm}
\begin{pf}
  a) By \Proposition{\ref{semi-cont-reg}} sets of the form $U_{k}=\{ x\
  \vert\ n_{inj}(x)<k\}$ are open.  By the noetherianness, this implies
  that there is an integer $N$ such that $k\geq N$ implies that
  $U_{k}=U_{N}$.  Since $V$ is separated at each point in $X$ this
  implies that $U_{N}=X$.  This proves the theorem when $X/k$ is
  regular.

  b) Assuming $X/k$ is not regular there exists a non-empty closed
  subset $X_s\subset X$ such that $\dim X_s < \dim X$ and $X_0 = X\setminus
X_s$ is
  regular \cite[III.10.5]{hartshorne}.  Provide $X_s$ with its reduced
  scheme structure; then $X_s/k$ again is a scheme of finite type. Let
  $j: X_s\to X$ be the associated closed immersion.  To prove the
  theorem we will use induction over $\dim X$.  By a),
  $n_{inj}^0:=\sup \{n_{inj}(x)\ \vert \ x\in X_0\}< \infty $.  If $\dim X =
  1$, its singular locus $X_{s}$ is a union of finitely many points
  $x_1, \dots , x_r$ so $n_{inj}= \max \{n_{inj}^0, n_{inj}(x_1), \dots
  , n_{inj}(x_r)\} < \infty$. We will extend this argument.

  Let $\xi $ be an associated point of $X_s\subset X$ and consider the
  following commutative diagram
  \begin{equation*}
      \xymatrix{
 V\ar[r]^\phi \ar[d] &  j^*(\Pc^n_{X/k})_\xi \ar[d]\\
k_\xi \otimes _k V \  \ar[r]^>>>>>{\dop_V^n(\xi)}  & \quad k_\xi
\otimes_{\Oc_\xi }\Pc^n_{X/k,\xi }(M),
}
  \end{equation*}
  where $\phi: V \to j^*(\Pc^n_{X/k})$ is the canonical map.  By
  \Definition{\ref{sep-def}}, there exists an integer $n_{inj}(\xi)$
  such that the map $\dop^n_V(\xi) : k_\xi \otimes_k V\to k_\xi \otimes_{\Oc_\xi
  }\Pc^n_{X/k, \xi }(M)$ is injective when $n\geq n_{inj}(\xi )$, hence $\phi$
  is injective when $n\geq n_{2}:=\max \{n_{inj}(\xi )\ \vert \ \xi $
  $\text{is an associated prime of } X_s\}$.


  This allows us to use the induction hypothesis.  We conclude that
  there is an integer $n_{1}$ such that the induced map $\gamma :
k_x\otimes_k V\to
  k_x \otimes_{\Oc_{X_s, x}} \Pc^{n_{1}}_{X_{s}/k,
    x}(j^{*}\Pc^{n_{2}}_{X/k}(M)) $ is injective for each $x\in X_s$.
  By \cite[Prop.  16.4.20, 16.7.9]{EGA4} we have a surjective map
  $$
  \alpha:j^*(\Pc^{n_{1}}_{X/k}(\Pc^{n_{2}}_{X/k}(M)))\to
  \Pc^{n_{1}}_{X_{s}/k}(j^{*}\Pc^{n_{2}}_{X/k}(M))
  $$
  and by \cite[Lemme 16.8.9.1, 16.7.9]{EGA4} a map
   \begin{displaymath}
\beta: \Pc^{n_{1}+{n_{2}}}_{X/k}(M)) \to
\Pc^{n_{1}}_{X/k}(\Pc^{n_{2}}_{X/k}(M)).
\end{displaymath}
Since, by functoriality, $\alpha(x)\circ \beta(x) \circ\dop^{n_{1}+n}_{X,
V}(x)=\gamma $ we
obtain that $\dop^{n_{1}+n_{2}}_{X, V}(x)$ is injective if $x\in X_s$.
Hence $\dop^{n}_{V}(x)$ is injective for each $x\in X$ if $n\geq
\max\{n_{inj}^0, n_{1}+n_{2}\}$.
    \end{pf}
    
    \begin{remark}\label{kam-remark}
      The assertion $n_{inj}(V)< \infty$ is proven in
      \cite{laksov-thorup} for smooth schemes $X/S$ with geometrically
      irreducible fibres.  In the language of
      \cite{kamran-Milson-Olver:invariant} the corresponding assertion
      is that $V$ is a ``regular'' subspace of $M$. Since they only
      consider varieties $X$ that are open subsets of some $\Rb^n$,
      which clearly are geometrically integral, arbitrary
      finite-dimensional subspaces $V$ of functions will be separable.
      In the analytic case $X_0$ may not be relatively quasi-compact
      and $n_{inj}(V)$ may thus be infinite for different reasons than
      in the algebraic situation.
  \end{remark}

\subsection{The surjectivity of the Taylor map}
Surjectivity properties of the Taylor map will play a role in
extending differential operators.

\begin{definition}\label{integer-def}
  Let $n_{surj}(x)$ be the largest integer such that $n\leq n_{surj}(x)$
  implies that $\dop^n_V(x)$ is surjective (so $\Cc^n_x =0$). The
  integer $n_{surj}(x)$ is the {\it jet order } \/of $V$ at $x$.
  Define $n_{surj}= \inf \{n_{surj}(x): x\in X\}$. Define also
  $n^1_{surj}$ to be the largest integer such that $\codim \supp \Cc^n
  \geq 2$ when $n\leq n^1_{surj}$.
\end{definition}
From (\ref{princ-fibre}) it is evident that if $x$ is a rational
point, then
\begin{displaymath}
n_{surj}(x) \leq \max \{n \ \vert \ \dim_k M_x/\mf_x^{n+1}M_x \leq
\dim_k V\}.
\end{displaymath}

The following well-known lemma illustrates the meaning of
$n_{surj}(x)\geq 1$. As mentioned in the introduction and as will be
described later every linear system $V\subset \Gamma(X,M)$ corresponds to a
rational map $X\to {\bf P}^n$, for some $n$. In the lemma the conditions
for this map to be a closed embedding is related to $n_{surj}$.

\begin{lemma}($k$ is alg.\ closed)\label{very-ample}
  Assume that $V\subset \Gamma(X,M)$ is a linear system on a
  non-singular projective variety $X$, where $M$ is an invertible
  sheaf. Then
   \begin{enumerate}
   \item $n_{surj}\geq 0 $ if and only if $V$ generates $M$;
   \item If $V$ defines a closed embedding e.g $V=\Gamma(X,M)$ and $M$ is
     very ample then for each point $x$ in $X$ we have $n_{surj}(x)\geq
     1$;
   \item If $n_{surj}(x)\geq 1$ and in addition $V$ separates points
     i.e.\ the canonical map $V\to k_x\otimes_{\Oc_x}M_x\oplus
k_y\otimes_{\Oc_y}M_y$ is
     surjective when $x$ and $y$ are different closed points in $ X$,
     then $V$ defines a closed embedding.
   \end{enumerate}
\end{lemma}

For the proof, consult \cite[Prop.II.7.3]{hartshorne}.

More generally $n_{surj}(x)\geq n$ is implied by the linear system being
$n$-jet ample (Cf.\cite{dirocco}). In the definition of the first
concept, in contrast to the second, only infinitesimal subschemes with
support in a point is considered. Hence it is in general weaker.

\begin{prop}\label{gen-weier} Let $M$ be a locally free
  $\Oc_X$-module. Then $n_{surj}\leq n^1_{surj} \leq N_{inj}$.
\end{prop}
\begin{pf} The first inequality is  obvious.
  Assuming $N_{inj}< n \leq n_{surj}^1$, then $\dop^n_V$ is an injective
  map that is surjective at points of height $1$ in the locus if
  non-singular points $X_0$. As $X_0$ is a regular variety
  $\Pc^n_{X_0/k}(M)$ is locally free so its depth at points of height
  $\geq 2$ is at least $2$; therefore $\dop^n_V$ is an isomorphism.  But
  since $n> N_{inj}$, the rank of $\Pc^n_{X_0/k}(M)$ is greater than
  the rank of $\Vc_{X_0}$, which gives a contradiction. Therefore
  $n_{surj}^1 \leq N_{inj}$.
\end{pf}

\section{Differential operators preserving $V$}

\subsection{Weierstrass points}

Applying $Hom_{\Oc_X}(\cdot , M)$ to the Taylor sequence
(\ref{w-sequence}) we get a map
\begin{equation}\label{dual}
  \Dc^n_X(M) \xrightarrow{W^n}   Hom_{\Oc_X}(\Vc_X, M) \cong
Hom_{k}(V,M).
\end{equation}
A (local) differential operator $P\in \Dc^n_X(M):=
Hom_{\Oc_X}(\Pc^n_{X/k}(M), M)$, induces a (local) map $\Hom_k(M,M)$,
or an action on $M$, by $m\mapsto P\cdot m:=P(d^{n}(m))$, where $d^{n}:M\to
\Pc^n_{X/k}(M)$ is the canonical map, described in 2.1. Then $W^n$
takes $P$ to the map $V\ni v\mapsto p\cdot v$, and we may describe the kernel of
$W^n$ as the annihilator $\Ann^n(V)= \{P \in \Dc^n_X(M)\ \vert \ P \cdot V =
0\}$. It is hence clear that $V$ defines an $\Oc_X$-coherent left
$\Dc_{X}(M)$-module
$$\Dc_{X}(M)/\Ann^n(V)\cong \bigcup_{n=0}^{\infty}W^n(\Dc^n(M))\subset
Hom_{\Oc_X}(\Vc_X,
M).$$
Note that its rank is less than or equal to $\dim V \cdot \rank M$.
We will first use this to give the upper bound on $N_{inj}(V)$,
referred to before, when $V$ is separated at each associated point.
This proposition is essentially proved in \cite{laksov-thorup}; we
have extended the result slightly and our setup perhaps gives a more
``conceptual'' proof.

\begin{prop}\label{v-1} ($X/k$ is a reduced scheme) 
  Assume that the map $\dop_V: \Vc_X \to \Pc_{X/k}^\infty$ is injective (see
  \Remark{\ref{sep-ass}}). Then $\dop_V^n$ is injective when $n\geq \dim
  V -1$, i.e. $N_{inj}(V)\leq \dim V -1$.  Also, the map $W^n_x$ is
  surjective at each associated point $x$ when $n \geq \dim V -1$.
\end{prop}

We need a well-known fact.
   \begin{lemma}\label{deg1-gen}
     If $K/k$ is a finitely generated field extension of
     characteristic $0$, then the ring of differential operators
     $\Dc_{K/k}$ is generated by its sub-space of first order
     differential operators.
   \end{lemma}
   By the references in the proof of \Lemma{\ref{locfree}} $K/k$ is
   differentially smooth (char. 0 ). Then the proof that $\Dc_{K/k}$
   is generated by $\Dc^1_{K/k}$ follows immediately from \cite[Th.
   16.11.2]{EGA4}.  (One may also notice that $\Dc_{K/k}(M)\cong
   \Dc_{K/k}\otimes_K End_K(M)$ is a matrix algebra of the ring $\Dc_{K/k}$,
   hence it is generated by matrices of differential operators in
   $\Dc^1_{K/k}$.)

\begin{pfof}{\Proposition{\ref{v-1}}}
  The map $\dop^n_V$ is injective if $\dop^n_{V,x}$ is injective for
  each associated point $x$, and then $\dop^n_{V,x}$ is injective if
  and only if $W^n_x$ is surjective since $\Oc_x$ is a field.  It
  suffices therefore to prove the following: Let $K/k$ be a field
  extension and $M$ a finite-dimensional $K$-linear space with a
  finite-dimensional $k$-linear sub-space $V\subset M$.  Then if the map
  $W:\Dc_{K/k}(M)\to Hom_k(V,M)$, $P\mapsto (v\mapsto P\cdot v)$ is
surjective, it
  follows that the map $\Dc^n_{K/k}(M)\to Hom_{k}(V,M)$ is surjective
  when $n\geq \dim_k V -1$. Clearly, it suffices to prove this when $M=K$
  and $V\subset K$. Moreover, if $W$ is surjective, since $\dim_ k V < \infty $
  there exists a subfield $K_1 \subset K$ that is finitely generated over
  $k$ such that $\Dc_{K_1/k}\to Hom_k(V, K_1)$ is surjective. Then if
  $\Dc_{K_1/k}^n\to Hom_k(V, K_1)$ is surjective it follows that
  $\Dc^n_{K/k}\to Hom_k(V, K)$ is surjective, since any element in
  $\Dc_{K_1/k}^n$ can be lifted to an element in $\Dc_{K/k}^n$. We can
  therefore also assume that $K$ is finitely generated over $k$.


  The space $Hom_k(V,K)$ is a module over $\Dc_{K/k}$, where a
  differential operator $P$ acts on $\phi: V \to K$ by $(P\cdot \phi)(v)=
  P(\phi(v))$. One then has $\Imo W^n = \Dc^n_{K/k}\Imo W^0 $ and $\dim_K
  \Imo W^0 = 1$.  Assume $\Imo W^n = \Imo W^{n+1}$, i.e. $
  \Dc^{n+1}_{K/k}\Imo W^0 = \Dc^n_{K/k}\Imo W^0 $.  Hence by
  \Lemma{\ref{deg1-gen}}
     \begin{displaymath}
        \Dc^n_{K/k}\Imo W^0 = \Dc^1_{K/k}\Dc^n_{K/k}\Imo W^0,
     \end{displaymath}
     and hence also that $ \Dc^n_{K/k}\Imo W^0 = \Dc_{K/k}\Imo W^0 $.
     By assumption $\Imo W = \Dc_{K/k}\Imo W^0 = Hom_k(V,K)$, and
     $\dim_K Hom_k(V,K) = \dim_k V < \infty$; hence there exists a smallest
     integer $n_0$ such that $\Imo W^n = Hom_k(V,K)$ when $n\geq n_0$.
     We have
     \begin{displaymath}
\Imo W^0 \subsetneq \Imo W^1 \subsetneq \cdots     \subsetneq \Imo
W^{n_0}= Hom_k(V,K),
   \end{displaymath}
   so the dimensions increase at each step.  Hence $n_0 +1$ $ \leq \dim_K
   Hom_k(V,K) = \dim _k V$, i.e. $n_0 \leq \dim_k V -1$.
%
%
\end{pfof}
A {\it gap}\/ for $V$ at a point $x$ in $X$ is an integer $i$ such
that $\rko \dop^i_V(x) > \rko \dop^{i-1}_V(x)$, and the {\it gap
  sequence}\/ of $V$ at $x$ is the set of gap integers.  We see from
the proof of \Proposition{\ref{v-1}} that the gap sequence at an
associated point $\xi $ of a separated sub-space $V\subset \Gamma(X,M)$ is
$1,2, \ldots
, n_{inj}(\xi)$ where $n_{inj}(\xi )\leq \dim V -1$.  (If $\trdeg_k K=1$,
then $N_{inj}=\dim V -1$.) Assuming $X$ is irreducible, one says that
a point $x$ in $X$ is a{\it Weierstrass point for $V$}\/ if its gap
sequence is different from the generic gap sequence. From
\Proposition{\ref{semi-cont-reg}} it follows that a point $x$ on a
non-singular variety $X/k$ is a Weierstrass point if and only if
$n_{inj}(x)> N_{inj}$.

For each integer $j\geq 0$ we have a short exact sequence
\begin{displaymath}
  0 \to \Vc \xrightarrow{\dop^{N_{inj}+j}_V}
\Pc^{N_{inj}+j}_{X/k}(M) \to \Cc^{N_{inj}+j}\to 0
\end{displaymath}
and can define the sub-sets
\begin{equation}\label{W-filtr}
W_j:=W_j(V)= \{x\in X \ \vert \  n_{inj}(x)> N_{inj} +j \};
\end{equation}
put also $W= W_0(V)$.  The closed set $W(V)$ (by the semi-continuity
of $n_{inj}(x)$) can be regarded as a set of {\it Weierstrass points
  on $X$ for $V$}, and $W_j(V)$ is its subset of {\it Weierstrass
  points of order j}.  We have
  \begin{displaymath}
   \emptyset = W_{n_{inj}-N_{inj}}\subseteq  \cdots  \subseteq
W_j\subseteq  \cdots \subseteq  W_1\subseteq  W.
  \end{displaymath}
  We can express these sets as supports of a coherent $\Oc_x$-module.
\begin{prop}\label{w-points}  ($X/k$  is a regular variety)
  We have
  \begin{displaymath}
    W_j = \supp Ext^1_{\Oc_X}(\Cc^{N_{inj}+j}, \Oc_X),
  \end{displaymath}
  so in particular $W_j$, $j=0,1,\dots , n_{inj} - N_{inj} +1$ are
  proper closed subsets.
\end{prop}
\begin{pf}
  $ \supp Ext^1_{\Oc_X}(\Cc^{N_{inj}+j}, \Oc_X)\subseteq W_j(V) $: Assume
$x\notin
  W_j(V)$.  By \Proposition{\ref{splitting}}, $(1)\Rightarrow (2)$, the Taylor
  sequence $(\ref{w-sequence})$, with $n= N_{inj}$ and
  $\Kc^{N_{inj}}=0$, is split exact since $\Pc_{X/k}^{N_{inj}+j}(M)$
  is locally free.  Therefore $\Cc^{N_{inj}+j}_x$ is free over
  $\Oc_x$; hence, $\Cc^{N_{inj}+j}$ being coherent,
  $Ext^1_{\Oc_X}(\Cc^{N_{inj}+j}, \Oc_X)_x =
  Ext^1_{\Oc_x}(\Cc^{N_{inj}+j}_x, \Oc_x)= 0$, so $x\notin \supp
  Ext^1_{\Oc_X}( \Cc^{N_{inj}+j}, \Oc_X)$.

  $W_k(V) \subseteq \supp Ext^1_{\Oc_X}(\Cc^{N_{inj}+j}, \Oc_X)$: If
$x\notin \supp
  Ext^1_{\Oc_X}( \Cc^{N_{inj}+j}, \Oc_X)$, then (\ref{w-sequence}),
  localised at $x$, is split exact; hence by
  \Proposition{\ref{splitting}}, $(2)\Rightarrow (1)$, $n_{inj}(x)= N_{inj}+j$;
  hence $x\notin W_j(V)$.

  That $\supp Ext^1_{\Oc_X}(\Cc^{N_{inj}+j}, \Oc_X)$ is a closed
  proper sub-set of $X$ is clear since $\Cc^{N_{inj}+j}$ is coherent.
\end{pf}
By \Proposition{\ref{w-points}} a natural scheme structure on $W_j(V)$
is given by the coherent ideal $\Ann_{\Oc_X} Ext^1_{\Oc_X}(
\Cc^{N_{inj}+j}, \Oc_X)$.

\begin{remark}
  Another way to say that $x\in W_j$ is that there exists a non-zero
  vector $v$ in $V$ such that $\Dc_{X/k}^{N_{inj}+j}(M)_x(v)\subset \mf_x
  M_x$. If $V$ is a complete linear system of affine dimension $r$ for
  a very ample invertible sheaf $M$, with respect to an embedding in
  projective space $X\to \Pb^{r-1}(V^*)=\Proj \So_k(V)$, this means that
  there exists a hyperplane in $\Pb^{r-1}$ that has contact with $X$
  of order $N_{inj}+j$ at $x$ (an osculating plane).
\end{remark}
Notice that the ideals $\Ann_{\Oc_X} Ext^1_{\Oc_X}(\Cc^{N_{inj}+j},
\Oc_X)$ define sub-schemes on any scheme $X/k$ locally of finite type
(making $\Cc^{N_{inj}+j}$ and hence the ideal coherent).  Thus $\bar
W_j :=\supp Ext^1_{\Oc_X}(\Cc^{N_{inj}+j}, \Oc_X)$ is one candidate
for Weierstrass sets of order $j$ for general schemes $X/k$, and also
for general $\Oc_X$-modules $M$, although we cannot expect $\bar W_j =
W_j$. Notice also that if $X/k$ is noetherian and reduced, $n_{inj} <
\infty $ \Th{\ref{n-fin-sing}}, so $W_j= \emptyset $ when $j\gg 1$, while it is not
certain that $\bar W_j$ need be decreasing in the singular case.
\begin{remark}
  In \cite{laksov-thorup} Weierstrass points on non-singular
  irreducible $X/k$ were defined using rank conditions on
  $\dop^i_V(x)$, giving rise to a sequence of closed subsets $Z(w_j)=
  \{x\in X \ \vert \ \rko \dop^i_V(x) < \rko \dop^i_V(\xi) \text{ for some
  } i = 0, \dots , j-1 \}$, where $\rko \dop^i_V(\xi )$ is the rank at
  the generic point $\xi$ in $X$. It is straightforward to see that $
  Z(w_j)$ equals $W(V)$ when $j\gg 1$, but for small $j$ the increasing
  filtration $\{Z(w_j)\}$ bears no natural connection to the decreasing
  filtration \{$W_k(V)$\} of $W(V)$.  In \cite{ogawa} Weierstrass points
  were defined differently.  Let $X/k$ be a projective regular variety
  and $h$ be the largest integer such that $\rank \Pc^n_{X/k}\leq \dim_k
  V$.  Ogawa defines the sets
  \begin{displaymath}
W^O_i = \left\{
  \begin{array}{lll}
& = \{x\in X \ \vert \ \dop^i_V(x) \text{ is    not surjective }\}
&\\
    &=  \{x\in X \ \vert \ n_{surj}(x) < i\}&\text {for } i \leq h,\\
&&\\
 &= \{x\in X \ \vert \   \dop^i_V(x) \text{ is not injective}\}&\\
    &= \{x\in X \ \vert \ n_{inj}(x) > i\}  &\text {for } i > h.
  \end{array}\right.
  \end{displaymath}
  Where in the second line we have expressed these sets using the
  functions $n_{inj}(x)$ and $n_{surj}(x)$.  Ogawa does not prove that
  $W_i^O = \emptyset $ for sufficiently high $i$; this however was proven in
  \cite{laksov-thorup}. By \Theorem{\ref{n-fin-sing}} we get that
  $W_i^O= 0$ for sufficiently high $i$ also when $X/k$ is any reduced
  scheme locally of finite type. As noticed in \cite{laksov-thorup}
  the sets $W^O_i$ need not be proper subsets of $X$.

  We are uncertain whether one can have $n_{inj}(x)< N_{inj}$ when $x$
  is a singular point. This cannot however occur at rational points on
  a curve $X/k$ when $M$ is an invertible sheaf, for then $n\geq
  n_{inj}(x)$ implies that $n\geq \dim V -1$, by (\ref{princ-fibre}), so
  $n\geq N_{inj}$ \Prop{\ref{v-1}}.
\end{remark}

\subsection{$M$ is simple as $\Dc_X(M)$-module when $X$ is regular}
Let $V$ be a separated $k$-subspace of $\Gamma(X,M)$, so in particular
there exists an integer $N_{inj}$ such that $\dop^n_V$ is injective
when $n\geq N_{inj}$.  Applying $Hom_{\Oc_X}(\cdot , M)$ to the Taylor
sequence (\ref{w-sequence}) we get an exact sequence
\begin{multline}\label{dual-wronski}
  0 \to Hom_{\Oc_X}(\Cc_n, M) \to \Dc^n(M) \xrightarrow{W^n}
  Hom_{\Oc_X}(\Vc_X, M) \to\\ \to Ext^1_{\Oc_X}(\Cc_n, M) \to
  Ext^1_{\Oc_X}(\Pc^n_{X/k}(M), M) \to
\end{multline}
We can as above identify $Hom_{\Oc_X}(\Cc_n, M)$ with the annihilator
$\Ann^n(V)= \{ P \in \Dc^n(M)\ \vert\ P \cdot V = 0\}$, and we will see in
\Proposition {\ref{abs-simple}} $(4)$ below that, conversely, if $M$ is
a simple $\Dc$-module, then $V$ can be recovered from $\Ann V$.

Clearly, the Taylor map $d^n_{V,x}$ is surjective if and only if the
evaluation $W^n_x$ is injective, hence the function $n_{inj}(x)$ may
be defined directly using the map $W^n$. In fact, on a regular variety
Weierstrass points may equally well be studied using differential
operators, as demonstrated in the following result:
\begin{thm}     \label{key-lemma}
  ($X/k$ is a regular variety) Let $M$ be a locally free
  $\Oc_X$-module, $\Dc_X(M)$ its ring of differential operators and
  $V$ a separated $k$-sub-space of $M_x$.  Consider the mapping of
  stalks
   \begin{displaymath}
     W^n_x: \Dc^n(M)_x \to Hom_{\Oc_x}(\Vc_x, M_x )= Hom_k(V, M_x).
   \end{displaymath}
   Then
\begin{displaymath}
      n_{inj}(x)= \min \{n \ \vert \ W^n_x \text{ is surjective.}\}.
\end{displaymath}
\end{thm}
\begin{pf} Put $s(x)= \min \{n \ \vert \ W^n_x \text{ is
    surjective.}\}$.

  $s(x)\leq n_{inj}(x)$: Let $n\geq n_{inj}(x)$. Since $x\mapsto n_{inj}(x)$ is
  upper semi-continuous \Prop{\ref{semi-cont-reg}} $n\geq N_{inj}$ so
  $\dop^n_{V,x}$ is injective. We may therefore use the sequence
  (\ref{dual-wronski}), so $W_x^n$ is surjective if $
  Ext^1_{\Oc_x}(\Cc_x^n, M_x) = 0 $.  Since $X$ is regular and $M$ is
  locally free it follows that $\Pc^n_{X/k,x}(M)$ is free. By
  \Proposition{\ref{splitting}}, $n\geq n_{inj}(x)$ implies that
  $\dop^n_{V,x}$ is split injective, so $\Cc^n_x$ is free; hence
  $Ext^1_{\Oc_x}(\Cc_x^n, M_x) = 0$.

  $s(x)\geq n_{inj}(x)$: Let $n\geq s(x)$.  We first prove that $s(x)\geq
  N_{inj}$.  Let $N^* = Hom_{\Oc_x}(N_x, \Oc_x)$ when $N_x$ is an
  $\Oc_x$-module.  Consider the Taylor sequence (\ref{w-sequence})
  localised at the point $x$, and let $\Ic_x$ be the image of
  $\dop^n_{V,x}$, so we have a surjection $\phi:\Vc_x \to \Ic_x$.  Now
  since $\Vc_x$ is torsion free it follows that $\Vc_x = \Ic_x$ if
  also the dual map $\phi^*: \Ic_x^* \to \Vc_x^*$ is surjective, which
  implies that $n\geq N_{inj}$ since $X$ is irreducible. To see that
  $\phi^*$ is surjective, we first note that by assumption the composed
  map $\Dc^n_{X/k,x}(M)\to Hom_{\Oc_x}(\Ic_x, M_x)\to Hom_{\Oc_x}(\Vc_x,
  M_x)$ is surjective. Hence, since $M_x$ is free, the composed map
  $\Pc^n_{X/k}(M)^* \to \Ic_x^* \xrightarrow{\phi^*}\Vc_x^*$ is surjective,
  hence $\phi^*$ is surjective.

  If now $n\geq s(x)\geq N_{inj}$, then, from (\ref{dual-wronski}), we have
  a split short exact sequence
   \begin{displaymath}
       0 \to Hom_{\Oc_x}(\Cc^n_x, M_x) \to \Dc^n(M)_x
\xrightarrow{W^n_x}   Hom_{\Oc_x}(\Vc_{x}, M_x) \to 0
   \end{displaymath}
   and since hence $M_x$ is free, we have a split exact sequence
   \begin{displaymath}
       0 \to Hom_{\Oc_x}(\Cc^n_x, \Oc_x) \to
Hom_{\Oc_x}(\Pc^n_{X/k,x}(M),\Oc_x) \to  Hom_{\Oc_x}(\Vc_{x},
\Oc_x) \to 0
   \end{displaymath}
   That $X/k$ is regular implies that the $\Oc_x$-module
   $\Pc^n_{X/k,x}(M)$ is locally free and hence reflexive, so upon
   applying $Hom_{\Oc_x}(\cdot , \Oc_x)$ to the previous exact sequence we
   get that the sequence
     \begin{displaymath}
     0 \to \Vc_{x} \xrightarrow{\dop_{V,x}^n} \Pc^n_{X/k,x}(M) \to
\Cc^n_x \to 0
   \end{displaymath}
   is split exact. Then \Proposition{\ref{splitting}} implies that $n\geq
   n_{inj}(x)$, whence with $(1)$ this implies that $s(x)=
   n_{inj}(x)$.
\end{pf}
\begin{prop}  \label{abs-simple}
  Let $X/k$ be a regular and geometrically integral variety and $x$ a
  point in $X$.
  \begin{enumerate}
  \item The map $W_x: \Dc_{X/k,x}(M)\to Hom_{k}(V,M_x)$ is surjective
    for each finite-dimensional $k$-sub-space $V\subset M_x$;
  \item $M_x$ is an absolutely simple $\Dc_{X/k,x}(M_x)$-module;
  \item $End_{\Dc_{X/k,x}(M)} M_x= k$;
  \item Let $L$ be a left ideal of $\Dc_{X/k}(M)$ such that
    $\Dc_{X/k}(M)/L$ is abstractedly isomorphic to $Hom_k(V,M)$ for
    some finite-dimensional $k$- sub-space $V$ (as $\Dc_X(M)$-modules,
    where $Hom_k(V,M)$ has the structure given by $Pf(v)=P(f(v))$, for
    $P\in\Dc_{X/k}(M),\ \ f\in Hom_k(V,M)$). Then there is a
    finite-dimensional vector space $V_{1}\subset M$, of the same dimension
    as $V$, such that $L= \Ann V_{1}$.  Define the map
  \begin{displaymath}
S: Hom_{\Dc_{X/k}(M)}(\Dc_{X/k}(M)/L, M)\to M, \ f\mapsto f(1 \mod  L).
\end{displaymath}
Then $S$ is injective and $V_1$ may be selected as the image of $S$.
\end{enumerate}
\end{prop}
\begin{pf}
  That $(1-3)$ are equivalent is proven in \cite[2.6.5]{dixmier}.
  $(1)$ follows from \Theorem{\ref{key-lemma}}.

  $(4)$: Let $v^*_i$, $i=1, \dots , n$ be a basis of the dual space
  $V^*$.  Let $\phi : \Dc_X(M)/L\cong Hom_k(V, M)$ be an isomorphism and
  $p_i: Hom_k(V, M)\cong V^* \otimes_k M \mapsto M$ be the projection
$v^*_j\otimes m \mapsto
  \delta_{ij}m$; $\Dc_X(M)$ acts trivially on the first factor in $V^*\otimes_k
  M$.  Then $p_i \circ \phi : \Dc_X(M)/L \to M$ is a surjection. Put $m_i =
  p_i\circ \phi (1 \mod L)$ and $V_{1}= km_1 + \cdots + km_n$.  Then
clearly, $L \subset
  \Ann V_{1}$, $V_{1}$ has the same dimension as $V$ and by $(1)$
  $\Dc_X(M)/\Ann V_{1}= Hom_k(V_{1},M) \cong \Dc_X(M)/L$.  Therefore $L=
  \Ann V_{1}$.  Finally,
   \begin{multline*}
     Hom_{\Dc_X(M)}(\Dc_X(M)/L, M) \cong Hom_{\Dc_X(M)}(V_{1}^*\otimes_kM , M
     )\\ \cong Hom_{\Dc_X(M)}(M,M)\otimes_k V_{1}\cong V_1
\end{multline*}
where the last isomorphism follows from $(3)$.
\end{pf}
\begin{remark}
  If $X$ is not regular it is well-known that $\Oc_X$ need not be a
  simple $\Dc_X= \Dc_X(\Oc_X)$-module. Notice that the surjectivity in
  $(1)$ in \Proposition{\ref{abs-simple}} for large $n$ follows from
  the density theorem, if $k$ is algebraically closed and $M_x$ is
  simple holonomic over $\Dc_X$.  Thus we have a counterpart of the
  surjectivity $(1)$ for any simple coherent $\Dc_X(M)$-module, not
  necessarily coherent over $\Oc_X$. The study of such $\Dc_X$-modules
  form the most interesting part of $D$-module theory, so it has some
  appeal extending the study to holonomic $\Dc_X$-modules.
\end{remark}

The next Theorem follows from \Theorem{\ref{key-lemma}}, together with
the identification of $\Ker W^n$ with $\Ann^n V$, discussed earlier.
\begin{thm}\label{surjective-seq}
  Let $X/k$ be a regular variety, $M$ a locally free $\Oc_X$-module,
  $\Dc^n_{X/k}(M)$ its $\Oc_X$-module of differential operators of
  order $n$, and $V$ a finite-dimensional separated $k$-sub-space of
  $\Gamma(X, M)$ \Defn{\ref{sep-def}}.  If $n\geq n_{inj}(V)$, then one has a
  locally split short exact sequence
   \begin{equation}\label{whole-split}
     0 \to \Ann^n V \to \Dc^n_{X/k}(M)\to V^* \otimes_k M \to 0.
   \end{equation}
\end{thm}

\subsection{The sheaf $\Dc^V(M)$}
We will now study the sheaf of rings $\Dc^V(M)$ of differential
operators that preserve a finite-dimensional vector space $V\subset
\Gamma(X,M)$.
Notice that
\begin{displaymath}
Hom_{\Oc_X}(\Vc_X, M) = Hom_k(V, M)= V^*\otimes_k M,
\end{displaymath}
where $V^*= Hom_k(V,k)$, and $End_k(V) \subset Hom_{\Oc_X}(V,M)$.
\begin{definition}
\begin{displaymath}
\Dc^{V,n}(M)= \{P \in \Dc_{X/k}^n(M) \ \vert \ P\cdot V \subset
V\}=(W^*_n)^{-1}(End_k(V))
\end{displaymath}
and $\Dc^V(M)= \cup_{n\geq 0}\Dc^{V,n}(M)$ (this is a sheaf of
$k$-algebras).
\end{definition}

Letting $n$ be an integer $\geq n_{inj}$, one gets the short exact
sequences (\ref{surjective-seq}).  One can push out
(\ref{dual-wronski}) to the short exact sequence
\begin{equation}
   \label{eq:sur-eq}
     0 \to \Ann^n V \to \Dc^{V,n}(M) \xrightarrow{W^n} End_{k}(V) \to
0.
\end{equation}

For obvious reasons this sequence is locally split exact, so let us
describe explicitly a splitting. Put $r= \dim V$ and let $\{\hat L^i\}$
be a basis of $End_{k}(V)$. Then select $L^i_x \in \Dc^{V,n_{inj}}_x(M)$
such that $W^n(L^i_x)= \hat L^i$ and define a local splitting $\phi_n :
\Dc^{V,n}_x(M)\to \Ann^n V $, $P_x\mapsto P_x - \sum_i \alpha_i L^i_x$ (the sum
contains $r^2$ terms), where the coefficients $\alpha_i\in k$ satisfy the
equation $\sum_i \alpha_i \hat L^i = W^n(P_x)$ in $End_k(V)$.  (Notice that
the same $\alpha_i= \alpha_i(W^n(P_x))$ works for all $x$ in affine subsets of
$X$; see the proof of \Theorem{\ref{global-dec}} below.) We thus need
to compute $r^2$ differential operators $L^i_x$ to define a split.
\begin{remark}\label{pres-order}
  The splits $\phi_n$ do not preserve the order filtrations of $\Ann^n V
  $ and $\Dc^{V,n}(M)$ when $n< n_{inj}$. But if $n\geq n_{inj}$ they do
  preserve the order filtration in the strong sense that $\phi_n$ induces
  an isomorphism $\Dc^{V,n}/\Dc^{V,n-1}\cong \Ann^n V/ \Ann^{n-1} V$.
  Cf.\/ also \cite{kamran-Milson-Olver:invariant}.
\end{remark}

We collect in the following theorems our general results on the
irreducibility of the action of $\Dc^V$ on $V$.  For $M= \Ac_{\Rb^m}$
the sheaf of real-valued analytic functions it was proven in [loc.\ 
cit., Th. 4.8], expressed in our terminology, that the map $W^n(\Omega)$ is
surjective for an {\it arbitrary} \/ open subset $\Omega$ of some $\Rb^m$
when $V$ is a subs-space of $\Ac_{\Rb^m}(\Omega)$ such that $n_{inj}(V)< \infty
$ and when $n\geq n_{inj}$. 
The following result allows locally free modules $M$ of rank $\geq 1$ and
affine algebraic varieties defined over arbitrary fields of
characteristic zero. The proof is at this point almost a formality.
\begin{thm} \label{global-dec}
  Let $X/k$ be a regular affine variety and $M$ a locally free
  $\Oc_X$-module. Let $V$ be a separated subspace of $\Gamma(X,M)$.
  Then if $n\geq n_{inj}$ we have a short exact sequence
   \begin{equation}\label{global-split}
  0 \to \Gamma (X,\Ann^n(M))\to \Gamma (X,\Dc^{n,V}(M))
\xrightarrow{W^n(X)} \End_k (V)\to 0
   \end{equation}
   and hence an exact sequence
\begin{equation}\label{global-split-all}
  0 \to \Gamma (X, \Ann (M)) \to \Gamma (X,\Dc^{V}(M))
\xrightarrow{W(X)} \End_k (V)\to 0.
\end{equation}
\end{thm}

\begin{pf}
  By \Proposition{\ref{surjective-seq}} the sequence
  (\ref{whole-split}) is locally split as vector spaces over $k$,
  hence if $X$ is affine, also globally split by Serre's vanishing
  cohomology theorem. Therefore the push-out (\ref{global-split}) is
  also split exact. We get (\ref{global-split-all}), since a split
  $End_k(V)\to \Dc^{n,V}(X)$ also gives a split $End_k (V)\to
  \Dc^V(X)$.
%
\end{pf}
This theorem implies that any geometrically integral real affine
variety (i.e. \/ there are no rational functions $f$ such that $f^2=
-1$) has the property that each map $V \to V$ is realised by a global
differential operator on $M$.  Hodge algebra, the main component in
the proof of \cite{kamran-Milson-Olver:invariant} of the corresponding
result for open subsets of some $\Rb^n$, is here replaced by the fact
that affine sets have vanishing higher cohomology.

By instead considering real points $X(\Rb)$ on a variety $X/\Rb$ (or
replacing $\Rb$ by any ordered field), and its sheaf $\Oc_{X(\Rb)}$ of
rational functions we get a stronger result since $X(\Rb)$ often is
affine, e.g.\/ when $X/\Rb$ is quasi-projective \cite[Th.
3.4.4]{bochnak}; $(X(\Rb), \Oc_{X(\Rb)})$ is an example of a real
algebraic variety [loc. cit].  A section of $\Oc_{X(\Rb)}$ is locally
a quotient $p/q$ of regular functions where $q$ does not vanish in any
real point. In the same way differential operators $\Dc_{X(\Rb)}$ have
the form $\sum_\alpha f_\alpha \partial^\alpha$ where $f_\alpha \in \Oc_{X(\Rb)} $ and
$\partial_{x_i}(x_j)= \delta_{ij}$ for some regular parameters $x_i$.

Let $M_\Rb$ and $\Dc_{X(\Rb)}(M_\Rb)$ be the localisation of $M$ and
$\Dc_X(M)$ on $X(\Rb)$.
\begin{thm}  
  Let $X/\Rb$ be a regular variety such that its assocated real
  algebraic variety $(X(\Rb), \Oc_{X(\Rb)})$ is affine. Then we have
  an exact sequence
  \begin{displaymath}
  0 \to \Gamma (X(\Rb), \Ann (M_\Rb)) \to \Gamma (X(\Rb),\Dc^{V}_{X(\Rb)}(M_\Rb))
\xrightarrow{W(X(\Rb))} \End_k (V)\to 0.
  \end{displaymath}
\end{thm}

\begin{pf} Let $x\in X(\Rb)\subset X$. 
  $V$ is separated at rational points \Cor{\ref{cor-sep}}; hence there
  exists an integer $n_{inj}(x)$ such that $W^n_x$ is surjective when
  $n\geq n_{inj}(x)$ \Th{\ref{key-lemma}}. Hence $W^n_x$ is surjective
  when $x$ belongs to a neighbourhood of $x$ in $X(\Rb)$. By
  quasi-compactness there exists an integer $n_{inj}$ such that
  $W^n_x$ is surjective for each $x\in X(\Rb)$ when $n\geq n_{inj}$.  We
  have $\Dc(M)_{X,x}\cong \Dc(M_\Rb)_{X(\Rb), x}$ when $x\in X(\Rb)$.  Hence
  the associated map
  \begin{displaymath}
W^n_\Rb: \Dc^n_{X(\Rb)/\Rb}(M_\Rb) \to   Hom_\Rb(V, M_\Rb)
\end{displaymath}
is surjective when $n\geq n_{inj}$. This implies that the global map
$W^n_\Rb(X(\Rb))$ is surjective by the real version of
\cite[Prop.II.5.6]{hartshorne} (it follows essentially for the reason
that the inverse image functor $i^*$ for the map $i:X(\bf{R})\to X$
takes injective $\Oc_X$-modules to flasque sheaves).
\end{pf}

It is not difficult to extend this line of argument to any
real-analytic manifold.

Similarly to [loc.\ cit, Th. 4.14] one can now estimate the order of
the generators of the annihilator of a separated subspace $V$.
\begin{cor}\label{ann-gen} Put
  \begin{displaymath}
A= \Gamma(X, \Oc_X),\quad  U=\Gamma(X, \Dc_X(M))\quad
\text{and}\quad  U^1=\Gamma(X, \Dc^1_X(M)).
\end{displaymath}
Assume that $U$ is generated, as $A$-algebra, by $U^1$.  Let $n$ be an
integer such that the sequences $(\ref{global-split})$ and
$(\ref{global-split-all})$ are split exact.  Then the left ideal $\Gamma
(X, \Ann V) \subset U$ is generated by the $A$-module $\Gamma (X,
\Ann^{n+1} V)$.
\end{cor}
Of course, the assumption on $U$ is always satisfied locally if $X/k$
is regular.
\begin{pf} Put $J^n= \Gamma(X, \Ann^n V)$ and $U^n= \Gamma(X,
  \Dc^n_{X/k}(M))$.  Use induction over $m$. If $m=n+1$, clearly
   \begin{displaymath}
J^m \subset U J^{n+1}.
\end{displaymath}
Assume $m> n+1$.  Let $\phi_m: U^m \to J^m$ , $m=1,2,\dots $, be the splits
described after the sequence (\ref{eq:sur-eq}). Since $U^1$ generates
$U$ it follows that if $P \in J^m \subset U^m $, $m> n$, then
\begin{displaymath}
   P = \sum P_{i}^{(1)} P^{(2)}_i
\end{displaymath}
where $P_i^{(1)} \in U^{m_1}$ with $m_1 < m$ and $P^{(2)}_i\in U^{n+1}$.
By \Remark{\ref{pres-order}} $P_i^{(2)}-\phi_{n+1}(P_i^{(2)})= P^{(2)}_i
\mod U^n$, hence $P_i^{(1)}P_i^{(2)} - P_i^{(1)} \phi (P_i^{(2)})$ is a
differential operator of order $\leq m-1$, i.e.
\begin{displaymath}
  P- \sum_iP_i^{(1)} \phi_{n+1} (P_i^{(2)}) \in J^{m-1},
\end{displaymath}
so $P \in J^{m-1} + U J^{n+1}$, and by induction
   \begin{displaymath}
    J^{m-1} \subset U J^{n+1},
   \end{displaymath}
   hence $P \in UJ^{n+1}$.
\end{pf}

In general, in a global situation when \Theorem{\ref{global-dec}} is
not applicable, it is hard to decide when $V$ is simple as module over
the global differential operators $\Gamma(X, \Dc^V(M))$.  Still there are
interesting cases when one can prove simplicity.  It is for example
true in certain cases both for toric varieties and homogeneous spaces,
as will be discussed later. Since in the case where the ground-field
is algebraically closed, this is clearly equivalent, by the density
theorem, to $W(X)$ being surjective we can describe the differential
operators on $X$ that preserve $V$, in an obvious but important way:

\begin{thm}\label{simple-dec}
  Let $M$ be a quasi-coherent $\Oc_X$-module and $V$ be a subspace of
  $V\subset \Gamma(X,M)$ that is simple as $\Gamma(X, \Dc^V_X(M))$-module.
Then
   \begin{displaymath}
\Dc^V_X(M) = \Gamma(X, \Dc^V_X(M))+ \mathrm{\Ann_{\Dc_X(M)}(V)}
\end{displaymath}
(equality of sheaves, where $\Gamma(X, \Dc^V_X(M))$ is regarded as the
constant sheaf of global sections of $\Dc^V_X(M)$)
\end{thm}

\section{Completions}\label{completions}

Finite-dimensional vector spaces of functions $V$ occur often as
linear systems associated to invertible sheaves on proper varieties.
Furthermore the differential operators that preserve these spaces,
considered as spaces of functions on an open subset, stem in important
examples, from differential operators on the proper variety.  In this
section we will explore conditions on $V$ that ensure that this
situation occurs. First we will recall the well-known method of
constructing a proper variety from a finite-dimensional subspace of a
$k$-algebra. To simplify the use we have given full references and
some arguments. The reader already familiar with this can read the
definition below, and then skip to section 4.2.
\subsection{Linear systems}

To fix our situation and keep track of our maps we define a category
$\Cc$ as follows.  An object $A$ in $\Cc$ consists of the datum
$(X,\Oc_X, V\xrightarrow{i} \Gamma(X,M),M)$ where $(X,\Oc_X)$ is a variety
over the field $k$, $V$ is a finite-dimensional vector space over $k$
and $i$ is a $k$-linear injective map to the global sections $\Gamma(X,M)$
of a locally free module $M$. Let $A=(X,\Oc_X, V\xrightarrow{i}
\Gamma(X,M),M)$ and $B= (X',\Oc_{X'}, V'\xrightarrow{i'} \Gamma(X,M'), M')$ be
objects in $\Cc$.  A morphism $J:A \to B $ is a morphism of
$k$-varieties $j : X \to X'$, an isomorphism of $\Oc_X$-modules $\psi :
j^*(M')\to M$, and an induced surjective map of
linear spaces $F=\Gamma(\psi): V' \to V$.

If $j:X\subset X'$ is the inclusion of an open subset, then the restriction
$B_{\vert X}:=(X, \Oc_{X},V\to \Gamma(X,j^{*}M),j^{*}M)$ belongs to $\Cc$,
and there is clearly a morphism $B_{\vert X}\to B$.

\begin{definition}\label{complete-def}
  An object $B=(X',\Oc_{X'}, V'\xrightarrow{i'} \Gamma(X,M'), M')$ in $\Cc$
  is a {\it completion}\/ of $A=(X,\Oc_X,i: V\to \Gamma(X,M),M)$ if there is
  a morphism $A\to B$ in $\Cc$ such that
   \begin{enumerate}
   \item $j: X\to X'$ is an open immersion into a projective variety
     $X'$;
   \item $i': V' \to \Gamma(X' , M')$ is an injection;
   \item $F: V' \to V$ is an isomorphism.
  \end{enumerate}
  If in addition $i'$ is an isomorphism, then the completion is {\it
    full}.
\end{definition}
Thus a completion is a simultaneous extension of $X$ to a projective
variety and a locally free extension $\bar M$ of $M$, with the
condition that $V\subset\Gamma (X', M')$, and it is full if $V\cong\Gamma
(X', M')$.  If
$j:X\subset X'$ is the inclusion of an open subset of a projective variety,
we see that $B_{\vert X}$ thus has the completion $B$, and it is a
full completion if $B= (X',\Oc_{X'}, \Gamma(X,M')\xrightarrow{Id} \Gamma(X,M'),
M')$.

First, we will use the standard algebraic geometric methods to
construct a canonical completion $C(A)$, as appropriate in our
situation, cf.  \cite{hartshorne, EGA2}.  We will only consider only
objects in $\Cc$ of the form $A=(X, R, V \subset R, R)$ where $X=\Spec R$ is
integral affine of finite type over $k$.

\begin{defn}\label{proj-def}
  Let $V\subset R$ be a finite-dimensional vector sub-space, and let
  $$
  \gamma :\\ \So [V][t]\to R[t],
  $$
  where $\So [V]$ is the symmetric algebra, be the graded
  homomorphism between rings graded by degree in $t$.  Let $B_V$ be
  the sub-algebra of $S[V][t]$ that is generated by $tV$ and denote by
  $A_V$ the image $\Imo \gamma(B_V)$.  Define
  \begin{displaymath}
X_{V}:={\Proj}A_{V}
\end{displaymath}
and let ${\mathcal L}_{V}= \Oc_{X_V}(1)$ be the associated invertible.
\end{defn}

Note that $X_{V}$ is a closed projective sub-variety of ${\Proj}(\So
[tV])={\bf P}^{(n-1)}$, if $\dim_k V=n$ and that ${\mathcal L}_{V}$ is
the pullback of $\Oc_{{\bf P}^{(n-1)}}(1)$, and hence very ample.

There is by definition an inclusion $i:A_{V}\to R[t] $, and hence a
rational map $i: X={\Proj} R[t]\to X_{V}$, defined on
$X-V_{+}(R[t]i((A_{V})_{+}))$ (Cf.\cite[3.1.7, 2.8.1]{EGA2}; $V_{+}()$
denotes the closed sub-variety defined by a graded ideal).  Now
\begin{displaymath}
R[t]i((A_{V})_{+})=R\oplus VRt\oplus VRt^2\oplus\ldots
\end{displaymath}
($VR$ is the ideal generated by $V$) and by \cite[2.3.13]{EGA2},
$V_{+}(R[t]i((A_{V})_{+}))$ is empty if and only if every element of
$R[t]_{+}$ has a power that is contained in $R[t]i((A_{V})_{+})$.
Applying this to $t\in R[t]_{+}$ we find that $VR=R$.  Hence $i$ is
defined on the whole of $X$ if and only if the following holds:

\noindent {\it Condition I}:\/
The ideal $VR=R$. (This is trivially true in the case that $V$
contains an unit.)

The map $i$ is furthermore an open immersion, if in addition the
following condition holds (it is not the most general condition
possible, but sufficient for our purposes):

\noindent {\it Condition II}:\/
\begin{equation}
   ((A_{V})_{tv})_{0}\cong \gamma(S[Vv^{-1}])= R_{v},
    \label{eq:VR}
\end{equation} for a set of $v=v_{i},\ i=1,\ldots,r$, such that
$(v_{1},\ldots,v_{r})=R$.  If $v\in V$ is a unit in $R$, and $\gamma$ is
surjective then actually $X=\Spec R\cong D_{+}(tv)$, since
$((A_{V})_{tv})_{0}=R$, and condition II is trivially satisfied.

That condition II implies that $i$ is an open immersion is easily
checked, since ${\Proj}A_{V}$ is constructed by glueing together the
affine spaces $D_{+}(tv):=\Spec ((A_{V})_{tv})_{0}\cong
\gamma(S[Vv^{-1}])$, for sections $tv\in tV\subset A_{V}$, and hence
the isomorphism in the condition only says that $i $ induces an
isomorphism $D_{+}(tv)\cong D(v)\subset \Spec R$, and that $D(v_{i}),\
i=1,\ldots,r$, cover $\Spec R$.

Assuming that conditions I and II hold, put
\begin{displaymath}
C(A):= ( X_V, \Oc_{X_V},V \xrightarrow{i}
\Gamma(X_V, \Lc_V), \Lc_V).
\end{displaymath}
There is a morphism $A\to C(A)$ given by $i :X\to X_{V},\
\psi:i^{*}(\Oc_{X_V}(1))\cong \Oc_{{ \Proj}R[t]}(1)\cong \Oc_{X}, F:V\
\xrightarrow{Id} V) $.

\begin{prop}\label{inverseV}
  Let $X=\Spec R$ be integral and assume that $A:=(X,\Oc_X,
  $ $ V\xrightarrow{i} \Gamma(X,\Oc_X),\Oc_X)\in\Cc$ satisfies condition I and
  II.
  \begin{enumerate}
  \item The canonical completion $C(A)$ is full if and only if the
    restriction map
    $$
    V=\Gamma({\bf P}^{(n-1)},\Oc_{{\bf P}^{(n-1)}}(1))\to \Gamma(X_{V},\Lc_{V})
    $$
    is surjective. This is true if e.g. $A_{V}$ is normal.
  \item There is a finite-dimensional vector space $V\subset
    V^{*}\subset R$ such that $A^{*}:= (X,\Oc_X, V^{*}\xrightarrow{i}
    \Gamma(X,\Oc_X),\Oc_X)$ has a full completion, in fact
    $(X_{V},\Oc_X,$ $
    \Gamma(X_{V},\Lc_{V})\xrightarrow{Id}\Gamma(X_{V},\Lc_{V}),\Lc_{V})$
    is such a full completion.
  \item Any full completion of $A$, $B= (\bar X, \Oc_{\bar X},
    \Gamma(\bar X, \Lc)\xrightarrow{Id} \Gamma(\bar X, \Lc) , \Lc)$,
    where $\Lc$ is very ample, satisfies $\bar X\cong X_V, \ \Lc\cong
    \Lc_V$.
      \end{enumerate}
\end{prop}

\begin{pf}    $(1)$:
  The restriction map may be identified with $\gamma$, which makes the
  first assertion obvious. Since $A_{V}$ is the homogeneous coordinate
  ring for the embedding of $X_{V}$ in ${\bf P}^{(n-1)}$, the
  assertion on normality is contained in \cite[Exc.II.5.15
  (d)]{hartshorne}.

  $(2)$: Put $V^{*}:=\Gamma(X_{V},\Lc_{V})_{\vert U}$.  Observe first that
  the restriction map
  $$
  \So^m(V)=\Gamma({\bf P}^{(n-1)},\Oc(m))\to \Gamma(X_{V},\Lc_{V}^m)
  $$
  is surjective for $m$ large enough, since $\Oc(1)$ is ample.
  Note further that $R[t]\cong \oplus_{i=1}^\infty \Gamma(X, \Lc^i)t^i$,
using the
  isomorphism $\psi : j^*(\Lc_{V})\to \Oc_X$, and that the restriction map
  hence gives an injection $\Gamma^{*}:= \oplus_{i=1}^\infty \Gamma (X_{V},
  j^{*}(\Lc_{V}^m)) \subset R[t]$, compatible with multiplication. The
  composition $\So^m(V)\to R[t]$ has the image $A_V$, which thus is the
  homogeneous coordinate ring of $X_{V}$, and coincides with the
  graded ring $\Gamma^{*}$ in high enough degrees.  Furthermore $A_{V^{*}}$
  is the image of the algebra $k[ V^{*}] $ generated by $V^{*}$ in
  $R[t]$, and it clearly contains $A_V$ and hence
  $(A_{V^{*}})_{m}=(A_{V})_{m}$, for $m$ large enough and consequently
  $\Proj(A_{V^{*}})\cong\Proj(A_{V})$, and $\Lc_{V}\cong\Lc_{V^{*}}$.

  $(3)$: This follows in a similar way. There is an injection
  $$
  \Gamma^{*}(\Lc):= \oplus_{i=1}^\infty \Gamma(\bar X, \Lc^i)t^i \subset R[t].
  $$
  Since $V=\Gamma(\bar X, \Lc)$, and $\Lc$ is very ample, the graded
  algebra $k[ V] $ generated by $V$ in $\Gamma_*$ coincides with $A_{V}$ in
  $R[t]$, and hence $X_{V}\cong \Proj(A_{V^{*}})\cong\Proj(\Gamma_*)\cong
\bar X$ (The
  last equality by very ampleness of $\Lc$).  Similarly we recover
  $\Lc$.
\end{pf}

\subsection{Extending differential operators}

We will give a condition on a finite-dimensional vector space $V\subset R$,
that ensures that a differential operator preserving $V$ actually
stems from a differential operator preserving an invertible module on
a proper variety.

\begin{thm} \label{ext-prop} Let $X$ be a variety, and $M$ a locally
  free $\Oc_X$-module.  Let $J: A \to B$ be a full completion, where
  $$A=(X, \Oc_{X}, i:V\to \Gamma (X, M),M)\ and \ B= (\bar X,\Oc_X,\bar i :
  \bar V\cong \Gamma(\bar X , \bar M), \bar M).$$
  Then the restriction map $\Gamma
  (\bar X, \Dc^n_{\bar X}(\bar M))\to \Gamma (X, \Dc^{n}(M)) $ takes its
  values in $\Gamma(X,\Dc_X^{n,V}(M))$, i.e.\/ we have an injective map
  \begin{displaymath}
    r_{\bar X,X}: \Gamma (\bar X, \Dc^n_{\bar X}(\bar M)\to\Gamma(X,
\Dc_X^{n,V}(M)).
  \end{displaymath}
  If $n\leq n^1_s$ \Defn{\ref{integer-def}}, i.e. $\codim_{\bar X} \supp
  \Cc^n_{\bar X} \geq 2 $, then $r_{\bar X,X}$ is an isomorphism.  If
  $\bar M$ is very ample, then in particular it is always an
  isomorphism for $n=1$.
\end{thm}

\begin{pf}
  We have $\Gamma (\bar X, \Dc_{\bar X}^n(\bar M)) = \Gamma (\bar X, \Dc_{\bar
    X}^{n, \bar V}) $. Since $\bar M$ is torsion free it follows that
  its ring of differential operators $\Dc_{\bar X}(\bar M)$ also is
  torsion free, so the restriction map $\Gamma (\bar X, \Dc_{\bar X}^n(\bar
  M)) \to \Gamma ( X, \Dc_{X}^n( M))$ is injective, and by the previous
  sentence its values preserve $\bar V \cong V$. This proves the first
  assertion.

  \noindent Assume $n\leq n^1_{surj}$.  An element $P$ in $
  \Dc^{n,V}(X)$ gives an element $\tilde P$ in $\Homo_{\Oc_{\bar
      X}}(\Vc_{\bar X},$ $ \bar M)=\Homo_k(V,\bar M)$.  The fact that
  $\tilde P$ comes from $P$ implies that generically the kernel $\Kc^n
  $ of $\dop^n_V:\Vc_{\bar X}\to \Pc^n_{\bar X/k}(\bar M)$ belongs to
  the kernel of $\tilde P $, but since $\Vc_{\bar X}$ and $\bar M$ are
  locally free, this gives $\Kc^n \subset \Ker \tilde P$.  Hence $\tilde P$
  gives an element $\bar P$ in $ \Homo_{\Oc_{\bar X}}(\Imo \dop^n_V
  ,\bar M)$.  Since $\dop^n_V$ is surjective at points of height $1$,
  $\bar P$ defines a map from $\Pc_{\bar X/k}^n(\bar M)$ to $\bar M$
  outside the codimension $\geq 2$ subset $\supp \Cc^n_{\bar X}$ of $\bar
  X$. Now $\bar M$ is locally free, so $\bar M_x$ has depth $\geq 2$ when
  $x\in \supp \Cc^n_{\bar X}$, hence $\bar P^n$ actually gives an
  element in $\Homo_{\Oc_{\bar X}}(\Pc^n_{\bar X/k}(\bar M),\bar M))=
  \Dc^n_{\bar X}(\bar M)(\bar X)$.

  The last assertion follows from the result earlier \Lem{\ref{very-ample}}
that very
  ampleness of a line-bundle implies that $n^1_{surj}\geq 1$ for the
  global sections.
\end{pf}

\begin{corollary}\label{van-anna}
  Keep the assumptions in \Proposition{\ref{ext-prop}} and assume also
  that $\dim_ k\Oc_X(X) = \infty$ (e.g. $X$ is affine). If $n\leq n^1_{surj}$,
  then $\Ann^n V= 0$.
\end{corollary}
\begin{pf}
  By \Proposition{\ref{ext-prop}} $\Gamma (X, \Dc_X^{n,V}(M))$ is
  finite-dimensional over $k$, hence $\dim_k \Ann^n V < \infty$.  However,
  $\Ann^n V$ is a torsion free $\Oc_X(X)$-module, implying that
  $\Ann^n V =0$.
\end{pf}

\section{Hidden Lie algebras}\label{hidden}

In this section we will consider the situation when $V\subset
A=\Oc_{X_0}$ is invariant under a reductive Lie algebra of
differential operators in $\Dc^1_{X_{0}}$.  General references for
representation theory and homogeneous spaces are
\cite{jantzen,kashiwara:rep, humphreys:algebraic-groups}.  We will
assume that $k$ is algebraically closed.
\subsection{An enlightening example}

Consider $X:={\mathbf A}^n$ and $V_{m}=\left< \prod_{i=1}^n x_{i}^{k_{i}}\
  \vert \ 0\leq k_{i} \ \text{and}\ \sum k_{i}\leq m \right>$.  We get that
conditions I and II of the preceding section are fulfilled and clearly
$\bar X_{V}={\mathbf P}^n$ and ${\mathcal L}_{V_{m}}={\mathcal O}
(m)$.  Thus
$$C(A)=({\mathbf P}^n,\Oc_{{\mathbf P}^n},V_{m}\cong\Gamma({\mathbf P}^n,\Oc
(m)), \Oc (m))
$$
is a full completion of
$$({\mathbf A}^n,\Oc_{{\mathbf A}^n},V\xrightarrow {} \Gamma({\mathbf
  A}^n,\Oc_{{\mathbf A}^n}), \Oc_{{\mathbf A}^n}).
$$

\begin{prop}\label{proj}
\begin{enumerate}
\item The Taylor map is an isomorphism
          \begin{displaymath}
\dop^m_{V}:\Oc_{{\mathbf P}^n}\otimes_k V_{m} \cong \Pc_{\Pb^n/k}^m(\Oc(m)).
        \end{displaymath}
      \item $ n_{inj}=n_{surj}=n^1_{surj}=m$.
           \end{enumerate}
        \end{prop}
         \begin{pf}
           By (\ref{princ-fibre}) the fibre of the Taylor map in $0\in
           {\mathbf A}^n$, with maximal ideal $\mf_0$, is $V_{m}\to
           \Oc_{{{\mathbf A}^n}}/\mf_0^{m+1}$, which is an
           isomorphism. $(2)$ now follows from this added to the fact
           that the Taylor map for homogeneous spaces and
           G-equivariant sheaves is equivariant, and hence surjective
           or injective in a fibre if and only if it has this property
           in every fibre. The equivariance of the Taylor map and the
           principal bundle is proved in
           \Proposition{\ref{G-linearized}}.
       \end{pf}

       Since $X=\mathbf{P}^n$ is a homogeneous space for the group
       $\Slo(n,k)$ there is a map $\beta: U(\Sl_{n})\to \Gamma(X,{\Dc(m)})$
(see
       \cite{kashiwara:rep}).

    \begin{cor}
      \cite{turbiner:bochner}
     \begin{enumerate}
     \item Any differential operator that has order less than $m$ and
       preserves the vector space $V_{m}$ of polynomials of degree
       less than $m$, is a polynomial in the differential operators
       $\partial_{x_{l}}, x_{k}\partial_{x_{l}},\ k,l=1,\ldots ,n$ together with
\\       $-\sum_{i=1}^{n}x_{i}x_{k}\partial_{x_{i}}+mx_{k},\ k=1,\ldots,n$.
     \item $\Dc^{V,1}\cong \Sl_{n}\subset \End_{k}(V)$
     \item $\beta(U(sl(n,k)))+\Ann V=\Dc _{V}$
   \end{enumerate}
\end{cor}

\begin{proof}
  Taking the global sections of $\Dc_X(m)$-modules $V= \Gamma(X, \Oc_X(m))$
  is a traditional way to construct a finite-dimensional simple
  $U(\Sl_{n})$-module.  Hence $V$ is also simple as
  $\Gamma(X,{\Dc(m)})$-module, so $W$ will be surjective, and
  \Proposition{\ref{simple-dec}} applies. This gives $(3)$.
  Furthermore $n^1_{surj}=m$, so \Corollary{\ref{van-anna}} is
  applicable and gives that $\Ann^n V=0$ if $n\leq m$, in particular for
  $n=1$; this implies $(2)$ (In fact it is easy to see that
  $\mathrm{Ann_{\D U}(V)}$ is the left ideal generated by the
  derivations $\partial^\alpha=\Pi \partial_{x_{i}}^{\alpha_{i}}$, where
$\alpha=(\alpha_{1},\ldots,\alpha)_{n}$,
  and $\sum_{1}^{n}\alpha_{i}\geq n+1$.) Finally, $(1)$ follows from the
explicit
  description of $\beta$; see \cite{musson,borel:Dmod, kashiwara:rep}.
\end{proof}

As a further easy corollary we can prove \Theorem{\ref{global-dec}} in
the case when $X= \Ab^n$ with a simple geometric argument. First a
lemma:

\begin{lemma}\label{simple-lemma}
  If $V$ is a subspace of $\bar V$ such that $\bar V$ is simple over
  $\Dc^{\bar V}$, then $V$ is simple over $\Dc^{\bar V, V} = \Dc^V\cap
  \Dc^{\bar V}$, the sub-sheaf of $\Dc^{\bar V}_X(M)$ consisting of
  sections that preserve $V$. Hence it is also simple over $\Dc^{V}$.
\end{lemma}
The idea of the easy proof is to choose a vector space complement
$\bar V=K\oplus V$, which induces an embedding $\End(V)\subset
\End(\bar V)$.

\begin{cor}
  Suppose that $V\subset \Oc_{\Ab^n}(\Ab^n)= k[x_1, \dots , k_n]$ is a
  non-zero finite-dimensional $k$-subspace.
   \begin{enumerate}
   \item $V \subset \Oc_{\Ab^n}(\Ab^n)$ is simple as $\Dc^V$-module.

   \item There is a completion $(\Pb^n, \Oc_{\Pb^n}, V_{m}\cong
     \Gamma({\Pb^n},\Oc (m)), \Oc_{\Pb^n} (m))$ for some positive integer
     $m$, of $(\Ab^n, \Oc_{\Ab^n}, V \subset \Oc_{\Ab^n}(\Ab^n),
     \Oc_{\Ab^n})$.  Furthermore
     $$\Gamma(\Pb^n,\Dc_{ \Pb^n/k}(\Oc(m))+\Ann V=\Dc^{V}.$$
   \end{enumerate}
     \end{cor}
\begin{proof}
  Embed $V$ in some $V_{m}=\{\sum_{0\leq \alpha_i ,\ \sum \alpha_i\leq
m}k_\alpha x^\alpha \ \vert \
  k_\alpha \in k \}$, where $\alpha = (\alpha_1 , \dots , \alpha_n)$,
$x^\alpha = x_1^{\alpha_1}\cdots
  x_n^{\alpha_n}$.  Then as above $\bar X_{V}={\mathbf P}^n$ and $\Lc
  _{V}=\Oc (m)$.  The enveloping algebra $\U(\Sl(n,k))$ gives
  naturally global differential operators on $\Oc (m)$ that makes
  $V_{m}$ a simple module.  Hence $V_{m}$ is simple also as a
  $\Dc^{V_m}$-module and \Lemma{\ref{simple-lemma}}, together with
  \Theorem{\ref{simple-dec}}, gives the result.
   \end{proof}

\subsection{Representations of Lie-algebras}
From now on assume that $k$ is algebraically closed.  We will now show
how these results extend to other homogeneous varieties and
representations (as also has been done by \cite{turbiner:bochner},
without using homogeneous varieties).

Consider first representations of a reductive connected and
simply\--con\-nected (not really necessary) group $G$.  Each irreducible
representation is constructed as $V=\Gamma(G/B,\Lc)$ for some ample line
bundle on the Borel variety $G/B$, where $\Lc=\Lc (\lambda)$ is associated
to an unique integral and dominant character $\lambda$ of the torus
$T\subset B$
(recall that dominant means that
$\langle\lambda,\alpha\rangle\geq 0$ for all roots $\alpha$
belonging to a basis of the root system consisting of positive roots).

The Borel variety $G/B$ contains an open cell $U$ that is a $B$-orbit,
which is isomorphic to some ${\mathbf A}^n$, and hence there is an
inclusion $V\subset \Oc_{G/B}(U)=k[x_{1},\ldots,x_{n}]=:R$, which is
completely specified
by the condition that it sends a primitive vector to $1\in R$. Actually it
is covered by
affine cells $gU$, all isomorphic to ${\mathbf A}^n$, so there are
many possible embeddings of $V$ in $R$.  If $\lambda-\rho$ is integral,
dominant, and regular(where $\rho$ is half the sum of the positive
roots), the line bundle is very ample, so the procedure
is invertible; we have $X_{V}\cong G/B$ and $\Lc_{V}\cong \Lc$ and hence also
$\Gamma(X_{V},\Lc_{V})=V$. If $\lambda$ is integral dominant, but not regular,
there are simple roots $\alpha$ such that $<\lambda,\check \alpha_{i} >=0,\
i\in I$, and
they define a parabolic group $P_{I}$. Furthermore, if $\pi : G/B \to G/P$
is the canonical projection, then $V=\Gamma(G/B,\pi^{*}(\Lc))\iso
V=\Gamma(G/P,\Lc)$, for a certain line-bundle $\Lc$. In this case
$X_{V}=G/P$, even if we start with $V\subset R$. Note how the proper variety
$G/B$ unites many different possible choices of open sub-varieties and
different vector spaces $V$.

Since $\frak{g} $ will act as 1st order differential operators (a
twist of a derivation) on $R=\Oc _{G/B}(U)$, and $V$ is irreducible as
$\frak{g} $-module, the situation of \Proposition{\ref{simple-dec}}
obtains. It is also well-known what the ring of global differential
operators of $\Lc (\lambda)$ is, see \cite{beilinson-bernstein}.  Hence the
following is immediate.

   \begin{prop}
     Suppose that $V \subset R$ is the inclusion of $V=\Gamma(G/B,\Lc
     (\lambda))$ that is specified by sending a $B$-primitive vector
     to $1$ in $R$.
     There is a canonical surjection $r: \U(\frak{g})\to \Gamma(G/B, \Dc
     _{\Lc (\lambda)})$, whose kernel is the ideal generated by a maximal
     ideal of the center $Z$ of $\U(\frak{g})$.  Hence
     $$\Dc_{R/k}^{V}=r(\U(\frak{g}))+\Ann V,
     $$
     and if $n\leq n^1_{surj}$
     $$\Dc_{R/k}^{n,V}=r(\U^n(\frak{g})).$$
     Moreover, if $X= \Spec R$,
     then $X_{V}=G/B$ and $\Lc _{V} =\Lc (\lambda)$
       \end{prop}
       This description of $\Dc^V_{R/k}$ is the main result of
       \cite{turbiner:bochner}.

\subsubsection*{G-equivariant sheaves}

Let $G$ be an arbitrary algebraic group, and $X$ a $G$-variety, with
the action given by $\mu:G \times X \to X$. Recall the following
induced
descent diagram:

\begin{displaymath}
  \xymatrix{
G\times G \times
X\ar@<3ex>[r]^>>>>>{d_0}\ar[r]^>>>>>{d_1}\ar@<-3ex>[r]^>>>>>{d_2} &
G\times X \ar@<3ex>[r]^>>>>>{d_0} \ar@<-3ex>[r]^>>>>>{d_1} &X
\ar[l]_>>>>>{s_0}
}
\end{displaymath}

\begin{alignat*}{2}
  d_0(g_1,x)&=g_1^{-1}x & \qquad \quad d_0(g_1,g_2,x)&= (g_2,
g_1^{-1}x)\\
  d_1(g_1,x)&= x&    \qquad\quad d_1(g_1,g_2,x)&= (g_1g_2,x)\\
  s_0(x)&=(e,x)& \qquad\quad d_2(g_1,g_2,x)&=(g_1,x)
\end{alignat*}
Note that all maps here are flat.  An $\Oc _{X}$-module $\Mcc$ is
called $G$-equivariant (see \cite{van-den-bergh, mumford}), if there
is an $\Oc _{G\times X}$-module isomorphism $\alpha:
d_{1}^{*}(\Mcc)\cong
d_{0}^{*}(\Mcc) $, such that the descent conditions
$d_{0}^{*}\alpha\circ
d_{2}^{*}\alpha=d_{1}^{*}\alpha $ and $s_{0}^{*}\alpha =id_{\Mcc}$
are true.

More properly $(M,\alpha)$ is really the object that should be called
equivariant, and such objects form, in an obvious way, an abelian
category.  However we will abuse notation and call $M$ equivariant. As
a module over itself $\Oc_{X}$ is $G$-equivariant, with structure map
$\alpha$ the composition of the canonical isomorphisms $\alpha:
d_{0}^{*}(\Oc_{X})\cong \Oc_{G\times X}\cong d_{1}^{*}(\Oc_{X}) $.
Below we will need that this map is an a homomorphism of sheaves of
rings.

The fact that the principal part bundle of a $G$-equivariant module is
$G$-equivariant is folklore; since we have been unable to find a
reference, we include for convenience a sketch of the simple proof.

        \begin{prop}
           \label{G-linearized}Let $X$ be a $G$-variety.
           If $\Mcc$ is $G$-equivariant, the $\Oc _{X}$-modules
           $\Pc^n_{X/k}(\Mcc)$ have compatible structures as
           $G$-equivariant modules. If $V\subset \Gamma(X,\Mcc)$ is
           invariant under the $G$-action induced from the equivariant
           structure, $\Oc _{X}\otimes_{k}V$ is canonically
           equivariant and the Taylor map (\ref{w-sequence}) is an
           equivariant homomorphism.
        \end{prop}

        \begin{corollary}
          Let $X$, $G$, $\Mcc$, and $V$ be as above. The functions
          $n_{inj}(x)$ and $n_{surj}(x)$ are constant along
          $G$-orbits.  In particular, if $X=G/H$ is a homogeneous
          space, then $n^1_{surj}=n_{surj}$, and $n_{inj}(x)=N_{inj}$,
          for all $x\in X$.
           \end{corollary}

\begin{proof}The diagonal map $\Delta: X \to X\times_k X$, is
  tautologically $G$-equivariant with respect to the diagonal action
  of $G$ on $X\times_kX$. Hence the canonical surjection
  $s:\Delta^{*}(\Oc_{X\times_k X})\to \Oc_{X}$, with kernel
  $I_{\Delta}$, is an equivariant map of $G$-equivariant sheaves. Call
  the isomorphism(the identity map after canonical identifications)
  $\alpha: d_{1}^{*}\Delta^{*}(\Oc_{X\times_k X})\cong
  d_{0}^{*}\Delta^{*}(\Oc_{X\times_k X})$.  Then $\alpha$ induces an
  equivariancy map $\alpha: \Ker d_{1}^{*}s= d_{1}^{*}(I_{\Delta})\to
  d_{0}^{*}(I_{\Delta})=\Ker d_{0}^{*}s$ (recall the flatness). It
  clearly suffices to see that
  $\alpha(d_{1}^{*}(I_{\Delta}^k))=d_{0}^{*}(I_{\Delta}^k)$ to obtain
  the equivariancy map $ d_{1}^{*}(\Pc^n_{X/k})\to
  d_{0}^{*}(\Pc^n_{X/k})$; this however is clear since $\alpha$, as
  noted above, is an algebra homomorphism, such that
  $\alpha(d_{1}^{*}(I_{\Delta}))=d_{0}^{*}(I_{\Delta})$. To check that
  this induces a structure of an equivariant module is now an easy
  exercise, as well as proving that the Taylor map is equivariant.
  (This proof more generally shows that the infinitesimal thickenings
  of an equivariant module with respect to a $G$-subvariety are
  equivariant.)

  The proof of the corollary is immediate from the further fact that
  the restriction of an equivariant module to an orbit is an
  equivariant module on the orbit, which is a homogeneous space $X=G/H$.
  For equivariant modules on homogeneous spaces there is an equivalence of
  categories between the category of coherent $G$-equivariant sheaves
  on $X$ and finite-dimensional
  $H$-representations(\cite{iversen,jantzen}).  In one direction this
  equivalence is given by taking the fibre $\Mcc_{H}$, with a
  canonical induced $H$-action, of the module $\Mcc$ in the
  $H$-invariant point. If $x=g^{-1}H\in X$, the stabilisator of $x$ is
  $H^g$, and the fibre $\Mcc_{x}$ is an $H^g$-representation, and then
  the different representations are related by $\Mcc_{x}\cong
^g\Mcc_{H}$. This is also true for morphisms, and shows that
the rank of an equivariant homomorphism is the same in each fibre over a
  closed point of the same orbit. This implies that the rank is the same
for each
  point of the orbit, since $k$ is algebraically closed.
  In particular, an equivariant module is locally free.
     \end{proof}

     In the case where $G$ is a reductive group it is now possible
     give a rough estimate of $n_{inj}(V)= N_{inj}$ and $n_{surj}(V)$
     when $V= \Gamma(G/B, \Lc(\lambda))$, since any invertible sheaf
$\Lc(\lambda)$ is
     $G$-equivariant on $G/B$ for a simply-connected group
     $G$(\cite{iversen}).

\begin{corollary}Let $w_{0}$ be the longest word in the Weyl
  group, $R$ the positive roots relative to $B$, and $\check{\alpha}$ the
  coroot associated to $\alpha\in R$. If $\lambda$ is a dominant weight and
  $\lambda-w_{0}(\lambda)= \sum_{\alpha\in R}k_{\alpha}\alpha$, denote the
minimal value of $\sum_{\alpha\in
    R}k_{\alpha}$ by $k(\lambda)$.  If $X=G/B$ and $V = \Gamma(X,\Lc)$, has
highest
  weight $\lambda$, then
  $$
  k(\lambda)\leq N_{inj}\leq \sum_{\alpha\in R}\langle \lambda,\check
\alpha\rangle.
  $$
  Furthermore $\lambda-w_{0}(\lambda)\geq \sum_{\alpha\in
R}k_{\alpha}\alpha$, if $\sum_{\alpha\in R}k_{\alpha}\leq
  n_{surj}$.
\end{corollary}
\begin{proof}
  Take $x=B\in G/B$. Then
  $$Hom_{\Oc_x}(\Pc^n_{X/k}(\Oc(\lambda))_{x},k)\cong \U^{n}(\frak{g})
v_{\lambda}\subset
  \U(\frak{g})\otimes_{\U(\frak{b})} kv_{\lambda},$$
  where the last module is the
  Verma module with highest weight vector $v_{\lambda}$, and
  $\U^{n}(\frak{g})$ denotes the elements in $\U(\frak{g})$ of order
  less than or equal to $n$ (\cite{beilinson-bernstein}).  Hence the
  dual of the fibre of the Taylor map at $x$, is the vector space
  homomorphism $\U^{n}(\frak{u}^{-}) v_{\lambda}\to V_{\lambda}$, induced
by the
  surjection of the Verma module to the irreducible module.  The
  Taylor map is injective if and only if its dual is surjective if and
  only if a non-zero element of lowest weight $w_{0}(\lambda)$ is contained
  in the image, implying the first inequality for $N_{inj}$.
  $F^{r}=\U^{r}(\frak{u^-}) v_{\lambda}$ is a filtration of $V_{\lambda}$, and
  $grV_{\lambda}$ is a $S(u^{-})$-module, and $N_{inj}$ may be expressed as
  the least $n$ such that $\mf^{n+1}\subset \Ann \gr V_{\lambda}$, where
$\mf$ is
  the ideal generated by $\frak{u}^{-}$. For the standard basis
  $x_{\alpha}$, the relations $x_{\alpha}^{n_{\alpha}+1}v_{\lambda}=0$,
where $n_{\alpha}=\langle
  \lambda,\check \alpha\rangle$ are true (\cite{humphreys:lie-algebras}),
implying the
  second inequality.  The inequality for $n_{surj}$ is immediate.
  \end{proof}
  
  For root systems of type $A_{1}$ and $A_{2}$ one has $k(\lambda)=N_{inj}$
  while the upper limit above for $N_{inj}$ is strictly larger in
  general for $A_{2}$. It would be interesting to know the precise
  value in terms of the weight $\lambda$.

\subsubsection*{Hidden symmetry}

Let us now consider the problem whether the above situation is in some
sense the only case. We have the following result, which might be
epistemologically interpreted as strengthening our general philosophy
that the construction $X_{V}$ is worthwhile to pursue since it (under
some conditions, of course) detects hidden geometry, in this case the
underlying homogeneous space.

\begin{prop}
     \label{hidden-symmetry}
     Given is a finite-dimensional $k$-subspace $V$ in an affine ring
     $R=\Oc_{X_{0}}(X_0)$, such that $X_{0}\subset X_{V}\subset
\mathbf{P}(V^*)$.
     Assume that $\frak{g}\subset \Dc_{X_{0}}^1(X_0)=\Dc^1_{R/k}$ is a
     reductive Lie algebra, and that $V$ is a representation.  Then
     the action of $\frak{g}$ on $V$ may be integrated to an action of
     an algebraic group $G$, whose associated Lie algebra is
     $\frak{g}$.  This action may be canonically extended to $X_{V}$,
     in such a way that $\Lc_{V}$ is an equivariant invertible sheaf.
     Furthermore, $\frak{g}\subset \Gamma(X_{V},\Dc _{X_{V}}(\Lc_{V}))$,
and $V\subset
     \Gamma(X_{V},\Lc_{V})$.  In the special case that $V$ is irreducible
     and $\frak{g}$ is locally transitive on $X_{0}$ and there is a
     point $x_{0}\in X_{0}$, such that the kernel of $\frak{g}\to \Dc
     _{X_0/k,x_0} $ is a parabolic sub-algebra, we have $X_{V}\cong G/P $
     for some parabolic subgroup $P\subset G$ and $\Lc_{V}=\Lc (\lambda)$
for some
     integral dominant weight $\lambda $.
\end{prop}
\begin{proof}

  The action of $\frak{g}$ on $V$ extends to a semi-simple homogeneous
  action on $\So[V][t]$, such that $\gamma:\So [V][t]\to R[t]$, is
  equivariant, where $\frak{g}$ acts through the inclusion $\frak{g}\subset
  \Dc^1_{R/k}$ on $R[t]$ \Defn{\ref{proj-def}}. Hence there is a
  compatible homogeneous action on $A_V=\Imo \gamma(\So [Vt])$, and this
  semi-simple action may be integrated to an action of a
  simply-connected algebraic group $G$ that has $\frak{g}$ as its Lie
  algebra (see\cite{humphreys:algebraic-groups}). Hence there is a
$G$-action on
  $X_{V}$, and it is easily checked that $\Lc_{V}$ is an equivariant
  invertible sheaf. From this follows that $\frak{g}\subset \Gamma(X_{V},\Dc
  _{X_{V}}(\Lc_{V}))$, and $V\subset \Gamma(X_{V},\Lc_{V})$.  This is the first
  part of the proposition. Using the additional assumptions in the
  second part, we get a map $g\mapsto gx_{0}$, $ G\to X_{V}$, where the kernel
  has to be a parabolic subgroup $P$ of $G$.  Hence we have a closed
  immersion $\phi:G/P\to X_{V}$. The local transitivity implies that there
  is an open subset $U\subset X_{0}$ that is in the image of $\phi$, and hence
  $\phi $ is an isomorphism. Since $V$ is a subset of the irreducible
  $G$-module $\Gamma(X_{V},\Lc_{V})$, it must equal the last module, and we
  are in the situation described in the beginning of the section.
\end{proof}

\section{Toric varieties}
We continue to assume that $k$ is algebraically closed. If $V\subset
k[x_{1},\ldots,x_{n}]$ is a finite-dimensional vector space generated
by monomials, $X_{V}$ will be a toric variety.  Differential operators
that preserve such $V$ have been considered in an affine situation in
\cite{post-turbiner,turbiner:sl2, turbiner:bochner,finkel-kamran},
without using toric varieties.  There is however much known about
rings of differential operators on toric varieties - for example
systematic procedures to calculate generators- by the work of Musson
\cite{musson} (see also \cite{jones,musson-bergh}), and one aim of
this section is to illustrate the use of this, as well as other toric
techniques.  We will also calculate the injectivity order for global
sections of line bundles on non-singular complete toric varieties; we
will in particular describe the locus of Weierstrass points for line
bundles on Hirzebruch surfaces.

\subsection{Completion of monomial vector spaces}

Recall the construction of projective toric varieties, cf. \cite{fulton} for
details. Let $M= \mathbf{Z}^n$ be a lattice. For
$m=(m_{1},\ldots,m_{n})\in M$ define
$$
x^{m}=x_{1}^{m_{1}}\ldots x_{n}^{m_{n}}\in
k[x_{1},x_{1}^{-1},.,\ldots,x_{n},x_{1}^{-1}]=k[T], \quad k[x]= k[x_1,
\dots , x_n].
$$
The torus in the nomenclature is $\Spec k[T]$.  If $V$ is the
vector space generated by monomials $x^m\ , m\in P$, where $P\subset
M$ is a convex polytope, $X_{V}$ \Sec{\ref{completions}} is usually
denoted $X(P)$ and the associated equivariant line bundle
$\mathcal{L}(P)$, cf. \cite[Section 1.5]{fulton}.  The construction of
$X(P)$ is explicitly described as follows. Let $m_{i}$ be a vertex of
$P$ and let $M_{i}\subset M$ be the semigroup generated by the set $\{
p-m_{i}\ \vert\ p\in P\} $ (or the elements of $M$ that lie in the
angle at $m_{i}$ bounded by the codimension $1$ faces that meet at
$m_{i}$).  Then define $k[M_{i}]\subset k[T] $ as the algebra
generated by $x^{m}$, $ m\in M_{i}$, and $U_{i}=\Spec k[M_{i}]$.
Furthermore define $L_{i}$ as $k[M_{i}]x^{m_{i}}\subset k[T]$, with an
obvious inclusion $V\to L_{i}$.  These local data $U_{i}$ and $L_{i}$
glue in a way that is uniquely determined by the given inclusions into
$k[T]$; this completes the construction.  We also have $V=\oplus_{m\in
  P} k x^m=H^0(X(P),\Lc (P))$ (\cite[3.4]{fulton}).  The line bundle
$\mathcal{L}(P)$ will be very ample if and only if the following
condition is satisfied:

\begin{description}
\item The polyhedron $P$ is the convex hull of the points $m_{i}$, and
  for each $i$ the semigroup generated by the set $ \{ p-m_{i} \ \vert\
  p\in P \}$ is saturated ([loc.cit.]).
\end{description}

As proven in \cite{musson} (see also \cite{jones,musson-bergh}), there
are sufficiently many differential operators on an invertible very
ample sheaf on a projective toric variety $X(P)$ to make $\Gamma(X,\Lc
(P))$ a simple $R=\Gamma(X,{\D{\mathcal L}(P)})$-module.  Hence by
\Proposition{\ref{simple-dec}} every differential operator preserving
$V$ may be decomposed as the sum of an annihilator and a global
differential operator on the line bundle.

\begin{prop}
    \label{toric-dec}Let $X=\Ab^n$. Assume that
    $V=\oplus_{m\in P\cap M} k x^m\subset k[x]$, where $P$ is a convex
    polytope satisfying the very ampleness condition above.  Then
  \begin{enumerate}
  \item $V=\Gamma(X_{V},\Lc_{V})$, and the restriction map
    $$
    \Gamma(X_{V},{\D{\Lc_{V}}})\to \Dc_{k[x]/k}^{V}/\Ann V
    $$
    is surjective;
  \item The object $(X(P),\Oc_{X(P)}, V\cong \Gamma(X(P),\Lc (P)),
    \Lc(P))$ is a full completion of $(X,\Oc_{X}, V\subset
    k[x],\Oc_{X})$ \Defn{\ref{complete-def}}.
  \end{enumerate}
  \end{prop}

  Musson describes the ring of global differential operators
  explicitly as a quotient of a ring of invariants of differential
  operators on an open subset of the affine space (see also
  \cite{jones} for a procedure to calculate generators). In general
  this gives smooth varieties where the global
  ring of differential operators is not generated by first order
  differential operators, hence not ``Lie-theoretic'' in the sense
  of \cite{kamran-Milson-Olver:invariant}.

\subsection{Calculation of $n_{inj}$ and $n_{surj}$ for toric varieties
in the smooth case}

\subsubsection*{Calculation of $n_{inj}(x)$}

The condition that the projective toric variety $X(P)/k$ be
non-singular is that the following {\it basis condition}\/ holds at
each vertex $m_i$ of the strongly convex polytope $P$ (Cf.
\cite{fulton}): There are exactly $n$ edges $E_{ij}$ meeting at $m_i$,
and if $e_{ij}$ is a minimal element of $M$ connecting $m_i$ with
another point in $E_{ij}$, then the set $\{e_{ij}\}_{j=1}^n$ is a
basis of $M$, which we call {\it the basis at the vertex $m_i$}. It
generates $M_{i}$.  Let $x$ be a point in $X(P)$ with residue field
$k_x$.  By equivariance the function $n_{inj}(x)$ is constant along
each orbit of the torus on $X(P)$.  Letting $\xi_1, \dots , \xi_b$ be
the generic points of the orbits it suffices therefore to compute the
numbers $n_{inj}(\xi_i)$.  If $\xi_i$ specialises to $\xi_j$, then
$P_j$ is a face of $P_i$ and by semi-continuity, $n_{inj}(\xi_j)\geq
n_{inj}(\xi_i)$. We put $N_{inj}(P) = n_{inj}(\xi)$ if $\xi$ is the
generic point of $T$ (or of $X(P)$). These orbits are in 1-1
correspondence with the faces $P_{r}$ of the polytope.  Suppose that
$m_{i}\in P_{r}$, and that $P_{r}$ contains precisely the basis
vectors $e_{ij}\ , j\in J_{r}$, in the basis at $m_{i}$. Then the
orbit corresponding to $P_{r}$ is an open dense subset of the
subvariety defined by the ideal $(x_{i},i\notin J_{r})$. It is
contained in $U_{i}=\Spec k[M_{i}]$.  Note that to the vertices of the
polytope there corresponds closed invariant points, and the closure of
each orbit contains at least one such point.  Denote by $H_{r}= \{m_i
+\sum_{j\in J_r} a_{ij}e_{ij}\ \vert \ a_{ij}\in \Nb \}$; then the
generic point $\xi_{r}$ of the orbit corresponding to $P_{r}$ has
$k_{\xi_r}=k(x^{e_{ij}}\ \vert \ e_{ij}\in H_r)$.  If $m= \sum_{j=1}^n
m_{ij}e_{ij}$, set $d(m)= \sum_{j\notin J_r}m_{ij}$. In the
proposition below we will consider parallel translates $H=H_{r}+m$ of
$H_r$ at a distance $d_H=d(m)$ from $H_{r}$.  The intersection $H \cap
P \subset P$ is a strictly convex polytope, parallel to the polytope
$P_r$, and we will relate $N_{inj}(H\cap P)$, the generic injectivity
order of the polytope $H\cap P$, to $n_{inj}(\xi_{r})$. We first note
that the generic injectivity order $N_{inj}(P)$ is intrinsic to $P$
and the lattice it generates.
\begin{lemma}
  Let $P$ be a polytope in $M$ and let $N$ be the sublatttice
  generated by $P$. The generic injectivity order of $V= \oplus_{m\in
    P}k x^m $ considered as a subsapce of $k[M]$, denoted $N_{inj}(P)$
  above, is the same as the injectivity order of $V$ considered as a
  subspace of $k[N]$.
\end{lemma}
We leave out the straightforward proof.

\begin{proposition} \label{inj-orbit}
  Assume that $P$ is a strongly convex polytope satisfying the basis
  condition, and that $X(P)/k$ and $\Lc(P)$ are the corresponding smooth toric
  variety and very ample line bundle. Let $\{\xi_i\}$ be the generic points
  of the orbits of the action by the torus.
 \begin{enumerate}
 \item \begin{multline}\label{orbit-n} n_{inj}(\xi_r) = \max\{ N_{inj}(H\cap
     P) + d_H \\ \ \vert \ H =m+H_{r}\text{ is a parallel translate of
     } H_r, H\cap P \neq \emptyset \};
\end{multline}
\item Let $P^v$ be the set of vertices in $P$ and $\xi_i$ be a closed
  orbit, i.e.\/ a closed point, with corresponding vertex $m_i$ in
  $P^v$, then
\begin{multline*}
  n_{inj}(\xi_i) = \max \{\ \sum_{j=1}^n m_{i,j}\ \vert \
  m\in P^{v} \}\\ \text{
    (the maximal distance from $m_i$ to any other vertex)};
\end{multline*}
\item
  \begin{multline*}
    n_{inj}(V) = \max \{ \ \max \{\ \sum_{j=1}^n m_{i,j}\ \vert \ m\in P^v\} \
    \vert \ m_i  \in P^v \}\\ \text{(the maximal length between vertices of
      $P$)}
\end{multline*}
\end{enumerate}
\end{proposition}
\begin{pf}
  $(1)$: Identifying $k_{\xi_r}\otimes_{\Oc_{\xi_r}}\Lc_{\xi_r} =
  k_{\xi_r}$ and $\Dc_{X(P)/k,{\xi_r}}(\Lc)= \Dc_{X/k,{\xi_r}}$ we
  want to determine the smallest integer $n_{inj}({\xi_r})$ such that
  $l \geq n_{inj}({\xi_r})$ implies that the evaluation map
\begin{displaymath}
W_x^l:k_{\xi_r} \otimes_{\Oc_x}\Dc_{X/k,{\xi_r}}^l \to Hom_k(\oplus_{m\in P}
kx^m,k_{\xi_r})
\end{displaymath}
is surjective \Th{\ref{key-lemma}}. Let $\phi_m \in \oplus_{m\in P}
Hom_k(kx^m,k_{\xi_r})$ be defined by $\phi_m(x^{m'})$ $=
\delta_{m,m'}\in k_{\xi_r}$.  That $n\geq n_{inj}({\xi_r})$ means that
$\phi_m\in \Imo W^n_{\xi_r}$ for each $m\in P$.  Let $P_{r}$ be the
face of $P$ corresponding to $\xi_r$, assume that the vertex $m_{i}\in
P_{r}$ and set $\{ x_{1},\ldots,x_{d}\}=\{x^{e_{ij} }\ |\ e_{ij}\in
J_{r}\}$, while $\{ x_{d+1}, \dots , x_n\}$ correspond to basis
vectors that are not parallel to $P_{r}$.  Let $\partial_{i}\in
\Dc_{k[T]/k}$ be defined by $\partial_{i}(x_i)= \delta_{ij}$.  Put
furthermore $H_r = m_i + \sum_{j\in J_r}\Zb e_{ij}$, so $P_{r}\subset
H_r$.  The residue field of the orbit is $k_F:= k(x_1, \dots ,
x_{d})$. Let $k[\partial^{(1)}] = k[\partial_{1}, \dots , \partial_d]$
and $k[\partial^{(2)}] = k[\partial_{d+1}, \dots , \partial_n]$.  For
$s\in P $ there is a unique decomposition $x^s=x^{m(s)}x^{f(s)}$,
where $f(s)\in H_{r}$, $m(s)\in\sum_{j\notin J_r}m_{ij}(s) e_{ij}$,
and the differential operator $\partial^{m(s)}:=\Pi_{j\notin J_r}
(\partial_{j})^{m_{i,j}(s)}\in k[\partial^{(2)}] = k[\partial_{d+1},
\dots , \partial_n]$ satisfies $\partial^{m(s)}(x^s)=Cx^f\in k_F$ for
some $0\neq C\in k$.  Note also that $\Imo W(k_F \otimes_k
k[\partial^{(1)}]\partial^m)(x^s)=0$ if $m\neq m(s)$, so that $\Imo
W(k_F \otimes_k k[\partial])(x^s)=k_F \otimes_k
k[\partial^{(1)}]\partial^{m(s)})(x^s)$.  Hence for any operator
$Q_{s}$ such that $W(Q_s)x^t=\delta_{t,s}$, there exists
$Q_{s}^1\subset k_F \otimes_k k[\partial^{(1)}] $ such that
$W(Q_s)=W(Q_{s}^1\partial^{m(s)})$ and
$W(Q_{s}^1)x^f=\delta_{f,f(s)}$; conversely, assuming the last
relation for some $Q_{s}^1$ we get that $
W(Q_{s}^1\partial^{m(s)})x^t=\delta_{t,s}$. If $o(Q)$ denotes the
order of a differential operator, and we assume $Q_s$ has minimal
order such that $W(Q_s)= \phi_s$, we have $o(Q_s)=o(Q_{s}^1)+m(s)$.
If $H=m+H_{r}$ is a parallel translate of $H_r$, then $s\in H\cap P$
iff $f(s)\in H\cap P - m$ and hence $\{W(Q_{s}^1)\ \vert \ s\in H\cap
P\}$ clearly is a $k_{F}$-basis of $\Hom(\sum_{p\in H\cap P -
  m}kx^{p},k_{F})$.  Hence $\max\{l\ \vert \ s\in H\cap P \ \text{and
} \phi_s\in \Imo W^l_{\xi_r}\}=N_{inj}(H\cap P - m)+d_{H}$, and since
$N_{inj}(Q -m)= N_{inj}(Q)$ for any polytope $Q$ in $M$, this
completes the proof of $(1)$.

$(2)$: This follows immediately from $(1)$ since $H_i$ is a
point.

$(3)$: This follows from $(2)$, since $n_{inj}(V)= \max \{n_{inj}(\xi_i)
\ \vert \ m_i \in P^v \}$.
\end{pf}

\Proposition{\ref{inj-orbit}} reduces the problem to the computation
of $N_{inj}= N_{inj}(P)$ ( $= N_{inj}(V)$) when $P$ is a strongly
convex polytope satisfying the basis condition.  Letting $K= k(x_1,
\dots , x_n)$ be the function field of $X(P)$ we can regard $V$ as a
subset of $K$. We have $\Dc_{X/k, \xi } = \Dc_{K/k} = K[\nabla_1, \dots ,
\nabla_n]$, where $\nabla_i$ are defined by $\nabla_i (x_j) =
\delta_{ij}x_i$, and we put
$A=k[\nabla_1, \dots , \nabla_n ]\subset K[\nabla_1 \dots , \nabla_n]$; put
also $A_l= k[\nabla]^l=
\{P\in A \ \vert \ \deg P \leq l\}$, where $\deg $ is the ordinary degree of
a polynomial.  Clearly, $A$ preserves each subspace $kx^m \subset V$, so $W$
restricts to a map $\bar W : A \to \End_k^d V$ (diagonal maps) and $W^l$
restricts to a map
\begin{displaymath}
\bar W^l:A_l\to \End^d_k(V) = \oplus_{m\in P} \Hom_k(kx^m, kx^m)\cong
k^{\vert P \vert }.
\end{displaymath}
If $I_l$ is the image of $\bar W^l$ it follows easily that $KI_l$ is
the image of $W^l$; hence $\rko \bar W^l = \rko W^l$.  Putting $J_P=
\Ker \bar W =\Ann V \cap A$ we get $\rko W^l = \rko \bar W^l = \dim A/(J_P
+ \mf^{l+1})$ where $\mf = (\nabla_1, \dots ,\nabla_n)$.  In particular,
\begin{displaymath}
N_{inj}(P)=  \min \{ l \ \vert \ \frac {A}{J_P+ \mf^{l+1}} =
\frac {k[\nabla]}{J_P+ \mf^{l+2}} \}.
\end{displaymath}
(Recall that $\rko W^0 < \rko W^1 < \cdots < \rko W^l \cdots <
\rko^{N_{inj}(P)}
= \vert P \vert $; see the proof of \Proposition{\ref{v-1}}). The
above expression for $\rko W^l$ was first given in \cite{perkinson}
using the explicit form of the matrix of the Taylor map $d^l_{V}$ in
the natural bases.

If $Z$ is a line in $M$ we let $|P\cap Z|$ be the number of points in $P\cap
Z$.  Put $d^g(P)= \max \{ | P\cap Z| \ \vert \ $Z$ \text{ a line in
  $M$}\}$. By \Proposition{\ref{v-1}} $N_{inj}(P)\leq \vert P\vert$, and
\Proposition{\ref{inj-orbit}} implies that actually $N_{inj}(P)$ is
bounded above by the maximal distance between the vertices of $P$
(computed in bases at the different vertices).
\begin{prop}\label{lem-N}
  Let $P$ be a subset of $M$. Then $ d^g(P) -1 \leq N_{inj}(P)$.
\end{prop}
\begin{pf}
  Let $Z$ be a line in $M$ such that $\vert Z\cap P\vert  =d^g(P)$.  Choose a
  point $m\in Z\cap P$ and let $P_m\in k[\nabla]$ satisfy $P_m x^{m'}=
\delta_{mm'}$.
  Then $P_m$ has a non-zero reduction $P^A_m$ in $A= k[\nabla]/I_Z$, where
  $I_Z$ is the ideal of $Z$, and $P^A_m$ vanishes at $d^g(P)-1$ points
  along $Z$; hence $\deg P_m \geq \deg P_m^A \geq d^g(P)-1$.  This implies
  $N_{inj}\geq d^g(P)-1$.
\end{pf}

\subsubsection*{Calculation of $n_{surj}$}
Let $P$ be a strongly convex polytope. The length $l_{ij}$ of an edge
$E_{ij}$, connecting the vertices $m_i$ and $m_j$ is the smallest
integer such that $m_j- m_i = l_{ij}e_{ij}$ (see above for the
definition of $e_{ij}$).  Denote by $s(P)$ the minimal length of an
edge ($P$ is $s(P)$-convex).

For a smooth and proper toric variety, Di Rocco \cite{dirocco} has
proved that $\Lc(P)$ is $s$-jet ample if and only if $s(P)\geq s$, for
$X(P)$ smooth and proper. This implies immediately that $s(P)\geq
n_{surj}$. Her method is to use Cox's homogeneous coordinate ring, to
describe the toric variety; this is not necessary in our simpler case.
We do not give a complete computation of the function $x\mapsto
n_{surj}(x)$, the surjectivity order at $x$, since we do not need it,
but it should be clear that it is locally constant on each orbit of
the torus, and may be estimated accordingly.

 \begin{prop}
    \label{toric-n}
    \begin{enumerate}
    \item Suppose that $P\subset M$ is a strongly convex polytope
      satisfying the basis condition.  Then $n_{surj}= s(P)$ and hence
      equals the degree of jet-ampleness.
    \item Let $P$ be a strictly convex polytope. Then \\ 
      $n^1_{surj}=\min \{n^1(F)\ \vert \ \text{F a codimension 1 face
        of }\ P\}$.

  \end{enumerate}
\end{prop}

\begin{pf}
$(1)$:  $X(P)$ is smooth if and only if any semigroup $M_{i}$ is generated
  by a basis of $M$. Fix $i$ and let
   $k[M_{i}]\cong k[x_{1},\ldots,x_{n}]=:k[x]$.
  The principal parts for $k[x]$ may be described as $d: k[x]\otimes_k
  V_{l}\iso \Pc^l_{\mathbf{A}^{n}/k}$, where $V_{l}$ is the vector
  space generated by all monomials in $x$ with total degree at most
  $l$, and the map is the Taylor map.
  Recalling that $L(P)x^{-m_{i}}=k[M_{i}]$, the $l$th Taylor map
  restricted to $U_{i}$ may be described in $U_{i}$, as
  $$
  k[M_{i}]\otimes Vx^{-m_{i}}\to \Pc^l_{U_{i}/k}.
  $$
  By convexity $P$ will contain all
  $\sum_{m=1}^na_{ij_{m}}e_{ij_{m}}+m_{i}$ where $\sum _{j}a_{ij_{m}}\leq
s=s(P)$. In
  particular $V_{s}\subset Vx^{-m_{i}}$, and by the description of
  $\Pc^s_{U_{i}/k}$ this implies that the $s$th Taylor map is
  surjective.  Conversely, assume that the $ s$th Taylor map is
  surjective. Tensoring with $k[x]/(x)$, we obtain that $Vx^{-m_{i}}\to
  k[x]/(x)^{s+1}\to k[x_{m}]/(x_{m})^s$ is surjective.  This gives that
  $\{1,x_{m},\ldots,x_{m}^s\}\subset Vx^{-m_{i}}$, which is, interpreted in
terms
  of $P$, precisely the condition that $s(P)=s$.
  
  $(2)$: Again the Taylor map is an equivariant homomorphism between
  $T$-homogeneous sheaves and so the support of the cokernel and
  kernel will be unions of closures of orbits under the torus $T$.
  Hence to check the surjectivity in codimension 2, it suffices to
  check surjectivity at the orbits of codimension 0 or 1 orbits.  In
  fact, it suffices to prove surjectivity for all orbits of
  codimension precisely $1$, since the support is closed and the only
  open orbit contains all codimension 1 orbits in its closure.
  \end{pf}

\begin{cor}
  Assume the conditions in \Proposition{\ref{toric-dec}}. If $m\leq s(P)$
  then the restriction map gives an isomorphism
  $$\Gamma (X_{V},\Dc^m (\Lc_{V}))\to \Dc_{k[x]/k}^{m,V}$$
    \end{cor}
 \begin{pf}
   This is immediate from \Corollary {\ref{ann-gen}}.
\end{pf}

\subsubsection*{Hirzebruch surfaces}
We exemplify with line bundles on Hirzebruch surfaces.  The polytope
is in this case determined by the finite-dimensional vector space of
polynomials
\begin{displaymath}
V_{kl}^{r}:= \left< x^{i}y^{j}\ \vert \quad\ 0\leq i+rj\leq k,\ 0\leq
j\leq
l\ \right>  ,
\end{displaymath}
where $r,k,l$ are non-zero integers and $r\geq 1$. We will restrict
ourselves to the truncated case $k-lr\geq 0$. As noted already in
\cite{gonzales-hurtubise-kamran-olver:quantification}, this vector
space is the restriction of the global sections of the equivariant
line-bundle $\Ok {\Sigma_{r}}(k,l)$ on the Hirzebruch surface $\Sigma_{r}$ to
the affine space $\mathbf{A}^{2}\subset \Sigma_{r}$.
The differential operators of order 1 that preserve the vector space
are described in [loc.cit], and in \cite{finkel-kamran} a graphic
method is given to calculate the higher order differential operators
that preserve $V$. This graphic method is just a use of the obvious
bigrading, and as such a special case of the much more powerful
methods of Jones/Musson. Even in the special case of Hirzebruch
surfaces, the methods of the latter authors give fuller information on
the whole ring of differential operators.

In particular, it is known that $\Dc^{V_{kl}^{r}}$ is not generated by
differential operators of order less than $1$. Let us see what the
preceding theory and the literature on toric varieties tells us.  The
vertices of the polytope are $V_1=(k,0), V_2=(0,0), V_3=(0,l),
V_4=(k-lr,l)$; denote the edge between the first two vertices by
$E_{1}$, between the second and third by $E_{2}$, and so on.

\begin{prop}Let $k-lr\geq 0$. Then
  \begin{enumerate}
  \item $$X=X_{V_{kl}^{r}}=\Sigma_{r}=\mathbf{P}(\Oc_{\mathbf{P}^1}\oplus
    \Oc_{\mathbf{P}^1}(r))$$
    and ${\Lc}_{V_{kl}^{r}}=\Oc_{\Sigma_r}(l,k)$;
\item $n_{surj}(X)=n^1_{surj}(X)=\Min \{l,k-lr \}$;
\item $N_{inj}= k$, $n_{inj}(V_1)= k$, $n_{inj}(V_2)= k+l$, $n_{inj}
  (V_3)= k+l$, $n_{inj}(V_4)= k$, $n_{inj}(E_1)= k$, $n_{inj}(E_2)=
  k+l$, $n_{inj}(E_3)=k$, $n_{inj}(E_4)=k$. In particular, $n_{inj} =
  k+l$;
\item The Weierstrass subsets (see (\ref{W-filtr})) are $W(V_{kl}^r) =
  W_l(V_{kl}^r) = \Pb^1$ where $\Pb^1$ is the closure of the orbit of
  the edge $E_2$;
\item The ring $\Dc^{V_{kl}^{r}}$ is given up to $\Ann V_{kl}^{r}$
  by the differential operators
  \begin{multline}
    R=\Gamma(\Sigma_{r},\Dc({\Lc}))   \\    =
k[\partial_{x},x^{j}\partial_{y},x\pi,
\partial_{x}^{j}y(\nabla_{y})\pi(\pi+1)\ldots(\pi+r-j-1),x\partial_{x},y\partial
_{y}\ | \ j=0,1\ldots,r].
\end{multline}
Here $\pi:=x\partial_{x}+ry\partial_{y}-k$, and $\nabla_{y}:=y\partial_{y}-l$.
\end{enumerate}
\end{prop}

\begin{pf}
  $(1)$: See \cite{dirocco,fulton}.

  $(2)$: The length of the edges are $k,l,k-lr, l$, and hence by
  \Proposition{\ref{toric-n}}, $n_{surj}(V_{kl}^{r})=\Min \{l,k-lr \}$.
  Next consider $ n^1_{surj}(V_{kl}^{r})$.  At the edge $E_{1}$, local
  coordinates are $x=x^{(1,0)}$ and $y=x^{(0,1)}$, and $p(F_{1})$ is
  defined by $x=1,y=0$. It is easy to see that $V_{kl}^{r}\to
  k[x,y]/(x-1,y)^{s+1}$ is surjective if and only if $s\leq l $, since in
  this case $y^s$ must be in the image.  In the same way
  $n^1(E_{2})=l$, $n^1(E_{3})=l$ and $n^1(E_{4})=k-lr$.  Hence
  $n_{1}=n_{surj}(V_{kl}^{r})=\Min \{l,k-lr \}$.

  $(3)$: Apply \Proposition{ \ref{inj-orbit}} to get the values along
  orbits of codimension $\geq 1$ (edges and points of P), noting that
  $N_{inj}(P)= |P|-1 = d^g(P)-1$ if $P$ belongs to a line in $M$.
  \Proposition{\ref{lem-N}} gives $N_{inj} \geq d^g(P)-1 = k$, and
  since $n_{inj}(V_1)=k$, by semi-continuity, we have $N_{inj}=k$.

  $(4)$: This follows from $(3)$.

  $(5)$: Since this is a toric situation, we know by Musson that $
  V:=V_{kl}^{r}$ is an irreducible module over the finitely generated
  and Noetherian algebra of differential operators $\Gamma (X_{V},\D {L})$.
  Hence \Theorem{\ref{simple-dec}} applies. Finally, an explicit
  description of a set of generators of global differential operators
  on $X$ is given in \cite{jones}. The restriction to $U_{1}=\Spec
  k[x,y]$ of the global differential operators on the structure sheaf
  $\Lc=\Oc_{\Sigma_{r}}$ are calculated to be
  \begin{multline}
    R=\Gamma(\Sigma_{r},\Dc({\Lc}))   \\
    =
k[\partial_{x},x^{j}\partial_{y},x\pi,\partial_{x}^{j}y(\nabla_{y})\pi(\pi+1)\ldots(\pi+r-j-1),x\partial_{x},y\partial
    _{y}\ |\ j=0,1\ldots,r].
\end{multline}
Here $\pi:=x\partial_{x}+ry\partial_{y}$ and $\nabla_{y}=y\partial_{y}$.  \
One may either repeat
these calculations for an arbitrary line-bundle --- in Jones and
Musson's framework this is an easy, if laborious exercise --- or one
may use \cite[Theorem 4.9]{jones}, to see that redefining
$\pi:=x\partial_{x}+ry\partial_{y}-k$, and also
$\nabla_{y}:=y\partial_{y}-l$ in the expression
$P(j)=\partial_{x}^{j}y(\nabla_{y})\pi(\pi+1)\ldots(\pi+r-j-1)$, will give
that the above
expression for the ring of differential operators is valid for the
line-bundle $\Lc=\Lc_{V_{kl}^{r}}$. This follows since after
redefinition the differential operators on the right hand side are
easily seen to act on $\Lc$, and the associated graded rings to the
filtration by differential operator order are equal.
\end{pf}

\bibliographystyle{siam} \bibliography{$HOME/.TeX/mindatabas}
\def\cprime{$'$}

\end{document}